\documentclass{imsart}

\RequirePackage[OT1]{fontenc}
\RequirePackage{amsthm,amsmath,amssymb}
\RequirePackage[numbers]{natbib}
\RequirePackage[colorlinks,citecolor=blue,urlcolor=blue]{hyperref}
\RequirePackage{mhequ,mathtools}
\RequirePackage{graphicx}
\RequirePackage{enumerate}
\usepackage[outdir=./]{epstopdf}
\usepackage[margin=1.25in]{geometry}
\arxiv{arXiv:1711.05724}

\startlocaldefs
\numberwithin{equation}{section}
\theoremstyle{plain}
\newtheorem{theorem}{Theorem}[section]

\newcommand{\be}{\begin{equs}}
\newcommand{\ee}{\end{equs}}
\newcommand{\bone}{\mathbf 1}

\DeclareMathOperator{\BF}{BF}
\renewcommand{\L}{\mathrm{L}}
\renewcommand{\d}{\mathrm d}
\DeclareMathOperator{\TV}{TV}
\renewcommand{\H}{\mathrm H}
\renewcommand{\P}{\mathbf P}
\newcommand{\E}{\mathbf E}
\newcommand{\bb}[1]{\mathbb{#1}}

\DeclareMathOperator{\Multi}{Multinomial}
\DeclareMathOperator{\Gam}{Gamma}

\newcommand{\Nt}{N}

\endlocaldefs

\begin{document}
\allowdisplaybreaks

\begin{frontmatter}
\title{Exact Limits of Inference in Coalescent Models}
\runtitle{Exact inferential limits for the coalescent}

\begin{aug}
\author{\fnms{James E.} 
\snm{Johndrow}\ead[label=e1]{johndrow@stanford.edu, juliapr@stanford.edu}}
\and
\author{\fnms{Julia A.} 
\snm{Palacios}
\ead[label=u1,url]{https://statistics.stanford.edu}}

\runauthor{J. Johndrow and J. Palacios}

\affiliation{Stanford University\thanksmark{m1}}

\address{390 Serra Mall\\
Stanford, California, 94305\\
\printead{e1}\\
\phantom{E-mail:johndrow@stanford.edu, juliapr@stanford.edu\ }}

\end{aug}

\begin{abstract}
Recovery of population size history from molecular sequence data is an important problem in population genetics.
Inference commonly relies on a coalescent model linking the population size history to 
genealogies. The high computational cost of estimating parameters from these models usually compels researchers to select a subset
of the available data or to rely on non-sufficient summary statistics for statistical inference.
We consider the problem of recovering the true population size history from
two possible alternatives on the basis of coalescent time data. We give exact expressions 
for the probability of selecting the correct alternative in a variety of biologically 
interesting cases as a function of the separation between the alternative size histories,
the number of individuals, loci, and the sampling times. The results are applied to human population history. 
This work has significant implications for optimal design when the inferential goal is to test scientific hypotheses about population size trajectories in coalescent models with and without recombination. 

\end{abstract}

%

\end{frontmatter}

\section{Introduction}

Estimation of historical effective population size trajectories from genetic data provides insight into how genetic diversity evolves over time. Availability of molecular sequence data from different organisms living today and from ancient DNA samples together with the development of evolutionary probabilistic models \citep{wakeley_coalescent_2008}, has enabled reconstruction of effective population size trajectories of human populations over the past 300,000 years \citep{Gattepaillegenetics.115.185058,Palaciosgenetics},  Ebola virus over the course of the 2014 epidemic in Sierra Leone \citep{EbolaBeast} and Egyptian hepatitis C virus for over 100 years \citep{HCVLast}. 

Until recently, inference of effective population size trajectories was limited by scarcity of molecular sequence data such as single nucleotide polymorphisms (SNPs) and microsatellites. The increase in the total amount of genetic data obtained at different time points from a large number of individuals over large genomic segments (loci), and the development of more realistic probabilistic models, has led to a situation in which computationally tractable statistical inference is only available from non-sufficient summary statistics such as the site frequency spectrum (SFS) \citep{Sainudiin2011}, from small numbers of samples, or from short genomic regions \citep{drummond2012bayesian, griffiths_sampling_1994, Kuhner:1995vw, Li:2011ez, StephensDonnelly2000}. \citet{GaoKeinan} give an extensive list of methods. 

Accurate detection of change points in the effective population size trajectory is of scientific relevance in many applications such as the timing and length of the human expansion out-of-Africa \citep{GaoKeinan}, and extinctions of large mammals at the end of the Pleistocene epoch often attributed to the depredations of humans \citep{shapiro_rise_2004}. Rather than studying the statistical properties of different estimators, we consider how increasing the amount of genetic data increases our ability to distinguish between alternative hypotheses about population history under different evolutionary models. We evaluate the ability to detect change points by asking what the lowest achievable error rate is for classification between alternative hypotheses about population history with different change points --the Bayes error rate. Calculation of Bayes error rates allows us to answer such questions without the need to rely on any particular estimator. We give exact expressions for the probability that the optimal procedure can identify the truth, given coalescence data. We study the effect of adding more sequences and more loci under the coalescent model with and without recombination and under different demographic scenarios. 

\section{Coalescent evolutionary models}

\begin{figure}[h]
\centering
\begin{tabular}{cc}
\includegraphics[width=0.8\textwidth]{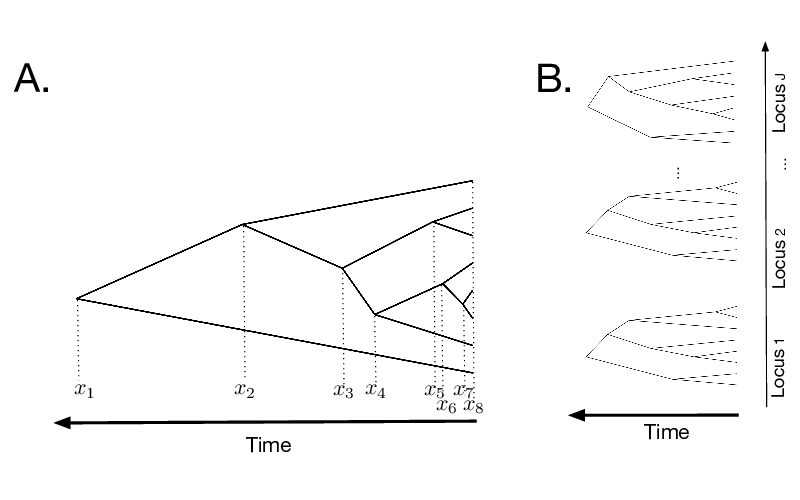} 
\end{tabular}
\caption{\textbf{A)} Genealogy of $n=8$ sampled individuals. $x_{i}$ is the time when two of $i+1$ extant lineages coalesce. \textbf{B)} Multiple genealogies along a chromosomal region.} \label{fig:coaltree}
\end{figure}

The standard n-coalescent \citep{Kingman:1982uj} is a generative model of molecular sequence of $n$ individuals sampled at random from a population of interest. In the single-locus neutral model, observed variation is the result of a stochastic process of mutations along the branches of the sample's genealogy; the genealogy is a timed bifurcating tree (Figure~\ref{fig:coaltree}A) that represents the ancestral relationships among samples. When moving back in time, two individuals find a common ancestor (coalesce) in the past with rate inversely proportional to the effective population size $N(t)$. Initially, the standard (homogeneous) n-coalescent assumed constant population size $N(t)=N$ and that sequences were sampled at the same time (t=0). Assuming a global mutation rate $\mu$, the parameter of interest is $\theta=2N\mu$. The standard neutral coalescent has been extended to variable population size $N(t)$ \citep{Slatkin:1991wx}, varying sampling times (heterochronous coalescent \citep{joseph_coalescent_1999}), and to account for population structure \citep{Beerli2001} and recombination \citep{griffiths1997ancestral}. 

Formally, the coalescent with variable effective population size $\Nt(t)$ \citep{Slatkin:1991wx} is an inhomogeneous Markov point process of $n-1$ coalescent times denoted by $x_{n-1},\ldots,x_{1}$. The process starts with $n$ individuals (lineages) at fixed time $x_{n}=0$  until  $x_{n-1}$ when  two of the $n$ lineages meet their most recent common ancestor. The process continues merging (coalescing) pairs of lineages until time $x_{1}$ when the remaining two lineages reach a common ancestor. The resulting realization is a genealogy with $n-1$ coalescent times like the one depicted in Figure \ref{fig:coaltree}A.
The conditional density of coalescent time $x_{k-1}$ 
is 
\be
f(x_{k -1}\mid x_{k},\Nt(t)) =\frac{C_k}{\Nt(x_{k-1})} \exp\left\{ - 
\int_{x_k}^{x_{k-1}} \frac{C_k}{\Nt(t)} dt \right\}
\ee
where $C_k = {k \choose 2}$ is the combinatorial factor depending on the number 
of possible ways that two lineages can coalesce given that there are $k$ 
lineages, and $\Nt(t)$ is the effective population size, a positive function of 
time. It follows that the complete likelihood is given by
\be \label{eq:IsoLikelihood}
L(x_1,\ldots,x_n \mid \Nt(t)) &= f(x_n) \prod_{k=n}^2 f(x_{k-1} \mid x_k, 
\Nt(t)) \\
&= \prod_{k=n}^2 \frac{C_k}{\Nt(x_{k-1})} \exp\left\{ - 
\int_{x_k}^{x_{k-1}} \frac{C_k}{\Nt(t)} dt \right\}
\ee
where again $x_n \equiv 0$ by definition.

In the coalescent model with recombination \citep{griffiths1997ancestral} looking backwards in time, lineages can either coalesce or recombine at a random position along the chromosome. When a lineage undergoes recombination, the lineage is split into two. The structure representing coalescent and recombination events is the ancestral selection graph (ARG). In the ARG, different chromosomal segments (loci) can have different genealogies and these genealogies are correlated (Figure \ref{fig:coaltree}B). \citet{McVean2005b} and \citet{MarjoramSMC} introduced Markovian approximations to the ARG called SMC and SMC' respectively. In the SMC, two genealogies at different segments separated by a recombination event are necessarily different, while in the SMC', these two genealogies are not necessarily different. Figure \ref{fig:coaltree}B shows an example realization of the SMC or SMC' process. In this paper, we analyze these approximations to the ARG from pairwise coalescent times. Derivation for larger sample sizes involve complicated likelihoods \citep{Palaciosgenetics} that are beyond the scope of this manuscript.

For $n=2$, let $x^{i}$ denote the pairwise coalescent time at the $i$-th locus; and let $J$ be the number of recombination events. In models with recombination, loci are contiguous chromosomal segments delineated by recombination, so we will also use $J$ to represent the number of loci. Under the SMC process, the transition density from $x^{i}$ to $x^{i+1}$, conditioned on a recombination event at locus $i+1$ is
\be\label{eq:smc}
f_{SMC}(x^{i+1}\mid x^{i}) = \frac{1}{x^{i}} \int^{x^{i+1}\wedge x^{i}}_{0}\frac{1}{ \Nt(x^{i+1})}q_{1}(u,x^{i+1})du,
\ee
where
\be
q_{k}(a,b)=\exp\left\lbrace-\int^{b}_{a}\frac{kdt}{\Nt(t)} \right\rbrace.
\ee
Given the current coalescent time $x^{i}$, a recombination breakpoint $u$ is sampled uniformly along the height of the tree $x^{i}$. At time $u$, one of the two branches is pruned with equal probability, and a new coalescent time $x^{i+1}$ is drawn with the standard coalescent rate.

Under the SMC' process, the transition density from $x^{i}$ to $x^{i+1}$, conditioned on a recombination event at locus $i+1$ is
\be \label{eq:smcprime}
f_{SMC'}(x^{i+1}\mid x^{i}) =\begin{cases}
\frac{1}{x^{i}}\int^{x^{i}}_{0}\int^{x^{i}}_{u} \frac{1}{\Nt(t)}q_{2}(u,t)
 dtdu & \left\lbrace x^{i+1}=x^{i} \right\rbrace \\
 \frac{1}{x^{i}\Nt(x^{i+1})}\int^{x^{i+1}}_{0} q_{2}(u,x^{i+1}) du & x^{i+1}<x^{i}\\
 \frac{1}{x^{i}\Nt(x^{i+1})}q_{1}(x^{i},x^{i+1}) \int^{x^{i}}_{0} q_{2}(u,x^{i}) du & x^{i+1}>x^{i}.
\end{cases}
\ee
Given the current coalescent time $x^{i}$, a recombination breakpoint $u$ is sampled uniformly along the height of the tree $x^{i}$. At time $u$, one of the two branches is selected with equal probability and split into two; one of the emanating branches follows the same trajectory back in time (old branch), while the other emanating branch can coalesce further back in time with any of the remaining two branches in the time interval $(u,x^i)$ with rate $2/\Nt(t)$. Conditional on failing to coalesce with any of the remaining branches in $[u,x^i]$, it coalesces with a branch emanating from the root 
at rate $1/\Nt(t)$ at some time $x^{i+1}>x^i$. The old branch is then removed. When the new branch coalesces back with the old branch, the resulting tree is the same as the original tree with coalescent time $x^{i}$. This event corresponds to the first case in equation \ref{eq:smcprime}.


The likelihood under SMC is given by
\be
\L(x^{1},\ldots,x^{J}\mid \Nt(t))=\frac{1}{\Nt(x^{1})}\exp \left\lbrace-\int^{x^{1}}_{0}\frac{dt}{\Nt(t)}\right\rbrace \prod^{J-1}_{i=1}f_{SMC}(x^{i+1}\mid x^{i}),
\ee
and the likelihood under SMC' for $n=2$ is obtained from the previous expression replacing $f_{SMC}$ by $f_{SMC'}$. 

\section{Bayes error rates in the standard coalescent}

We will start with the simple hypothesis testing setting
in which the two hypotheses are:
\be \label{eq:hypo1}
H_1 : \Nt(t) = a N_0, \quad T\le t \le T+S \\
H_2 : \Nt(t) = b N_0, \quad T \le t \le T+S
\ee
with $\Nt(t)$ equal under $H_1$ and $H_2$ outside the interval $[T,T+S]$; $a, b$ and $S$ are positive constants and $T\geq 0$ (Figure \ref{fig:OoAnJ}). Our
goal is to determine which hypothesis represents the true state of nature under which
the data were generated. For simplicity of notation, we associate the state of nature with
a parameter $\vartheta \in \{1,2\}$ such that $H_1 : \vartheta = 1$ and $H_2 : \vartheta=2$.
A (binary) Bayes classifier or decision rule $\vartheta(x)$ has
the form
\be
\vartheta(x) = \begin{cases} 1 & \BF_{12}(x) > 1 \\ 2 & \BF_{12}(x) < 1 \\ \xi & \BF_{12}(x) = 1 \end{cases} 
\ee
where $\BF_{12}(x)$ is the Bayes factor for $H_1$ vs $H_2$,  $x$ is an 
observation of a random variable $X$, and $\xi \sim \text{Bernoulli}(1/2)+1$. In the sequel, we drop the explicit argument and simply write $\BF_{12}$ in place of 
$\BF_{12}(x)$.
Thus, if $\vartheta(x)$ returns 1, we infer that the data were generated under $H_1$,
whereas if $\vartheta(x)$ returns 2, we infer that the data were generated under $H_2$.
In the case where 
each hypothesis is assigned prior 
probability of one half, the Bayes factor is exactly the likelihood ratio, and the probability of selecting $H_1$ is the probability that $\BF_{12} > 1$ plus half the probability that $\BF_{12}=1$. In this case, the probability of correct classification is
\begin{eqnarray*}
\P[\vartheta(X) = \vartheta]=& \frac{1}{2}\left[ \P(\log \BF_{12}>0 \mid H_{1})+\frac{1}{2}\P(\log \BF_{12}=0 \mid H_{1})\right] \\
& +\frac{1}{2}\left[\P(\log \BF_{12}<0 \mid H_{2})+\frac{1}{2}\P(\log \BF_{12}=0\mid H_{2})\right]
\end{eqnarray*}
When the prior is correct, the Bayes classifier is the optimal
classifier, so that the probability of correct classification using the Bayes classifier is the maximum achievable probability. The Bayes error rate is $1-\P[\vartheta(X) = \vartheta]$. As such, by studying the properties of the Bayes classifier, we obtain general limitations on inference for any classifier or test. 

We first define some notation. Let $X = (X_1,X_2,
\ldots,X_{n-1})$ be the random vector of coalescent times with distribution given by \eqref{eq:IsoLikelihood}.
When multiple genealogies are available, we will denote the random variable corresponding to the collection of all $J$ genealogies by $X^J$. Throughout, we abuse notation by writing $\P[\vartheta(X) = \vartheta]$ -- or $\P[\vartheta(X^J) = \vartheta]$ when $J>1$ -- to represent the
probability of correctly identifying the true state of nature.


The following theorems provide exact expressions for the probability of distinguishing between two hypotheses of the form \ref{eq:hypo1} from pairwise coalescent data under the coalescent with variable population size \eqref{eq:IsoLikelihood}.


\begin{theorem} \label{thm:BayesError1Locus}
Consider the simple hypothesis testing problem of the form \ref{eq:hypo1} when a single pairwise coalescent time is observed ($n=2$) and assign equal prior probabilities to both hypotheses. Then the probability of correctly distinguishing between the two hypotheses is:  
\be
\P[\vartheta(X) = \vartheta] = \frac12 + \frac12 e^{-\Lambda(T)} \left( e^{-\frac{ \delta \wedge S}{(a \vee b)N_0}} - e^{- \frac{\delta \wedge S}{(a \wedge b) N_0}} \right)
\ee
where
\be
\delta \equiv \frac{ab N_0}{b-a} \log \frac{b}{a} = \frac{ab N_0}{a-b} \log \frac{a}{b} \ge 0,
\ee
and
\be
\Lambda(T) \equiv \int_0^T \frac{1}{\Nt(t)} dt.
\ee

\end{theorem}

Proofs of all results 
can be found in the Appendix. 
Type I and type II error rates can be obtained from the conditional probability expressions derived in the proof, and by modifying our proof to consider a classifier that thresholds $\BF_{12}(x)$ at $\zeta(\alpha)$ for which $\P[\BF_{12}(X) > \zeta(\alpha) \mid H_1] = 1-\alpha$, one can perform power calculations for testing at level $\alpha$ where $H_1$ is designated as the null. We mention this as an obvious extension of our results, but do not pursue power calculations in the current work.

We now extend our previous result to the case when $J$ independent genealogies from \eqref{eq:IsoLikelihood} are available. When multiple loci or chromosomal segments are either coming from different chromosomes or from the same chromosome at distant locations, genealogies at those locations can be assumed to be independent. When $n=2$ and $J$ independent genealogies with likelihood (\ref{eq:IsoLikelihood}) are available, the sample configuration $L=(L_{1},L_{2},L_{3})$ of the $J=L_{1}+L_{2}+L_{3}$ pairwise coalescent times is $L\sim \text{Multinomial}(J,\mathbf{p}=(p_{1},p_{2},p_{3}))$, where $L_{1}$ is the number of pairwise coalescent times that fall in the interval $(0,T)$, $L_{2}$ is the number of pairwise coalescent times that fall in the interval $[T,T+S]$, $L_{3}$ is the number of pairwise coalescent times that are greater than $T+S$, and
\be
p_1 &= \P[X \le T] = 1-e^{-\Lambda(T)} \\
p_2 &= \P[T < X \le T+S ] = e^{-\Lambda(T)}-e^{-\Lambda(T+S)}
\\
p_3 &= \P[X > T+S] = e^{-\Lambda(T+S)}.
\ee
%
%
%
%
For this setting we have the following result
\begin{theorem} \label{thm:BayesErrorJLoci} Consider the simple hypothesis testing problem of the form \eqref{eq:hypo1} when $J$ independent pairwise coalescent times are observed ($n=2$, $J \geq 1$). The probability of correctly distinguishing between the two hypotheses is  
\be
\P[\vartheta(X^J) = \vartheta] &= \frac12 \P(L_2=0 \mid H_1) + \frac12 \sum_{(\ell_2,\ell_3) : \ell_2 > 0} \P(L=\ell \mid H_1) \P[W^*(\ell_2) > \ell_2 \delta - \ell_3 S \mid H_1, L=\ell] \\
&+ \frac12 \sum_{(\ell_2,\ell_3) : \ell_2 > 0} \P(L=\ell \mid H_2) \P[W^*(\ell_2) < \ell_2 \delta - \ell_3 S \mid H_2, L=\ell] \\
\ee
where $W^*(\ell_2) = \sum_{j=1}^{\ell_2} X^{j}_{*}$ is the sum of $\ell_2$ independent truncated coalescent times $X_*^j \in [0,S]$, each exponentially distributed with rate $(aN_0)^{-1}$ under $H_1$, and rate $(bN_0)^{-1}$ under $H_2$; $\delta$ is defined as in Theorem \ref{thm:BayesError1Locus}, and


\be
\ell \in \big\{\ell = (\ell_1,\ell_2,\ell_3) : \ell_j \in \bb N,\, \sum_j \ell_j = J \big\}.
\ee
is an element of the support of $\Multi(J,\mathbf{p}=(p_{1},p_{2},p_{3}))$. 
\end{theorem}

To obtain numerical results, we approximate the distribution function $\P[W^*(\ell_2)<t]$ by Monte Carlo. 
In the next section we apply these results to the problem of distinguishing between two hypotheses about the human expansion out-of-Africa.



\section{Human expansion out-of-Africa} \label{sec:OoA1}
Many of the statistical methods proposed over the last 15 years to infer effective population sizes from genetic data have been applied to human whole genomes \citep{Li:2011ez,MSMC,diCal,Palaciosgenetics,TerhorstSFS}. Several studies agree that non-African populations have experienced two severe bottlenecks, one attributed to the expansion out-of-Africa and the other attributed to the separation of Asian and European populations. There is, however, disagreement in the timing and length of such events. 

\begin{figure}[h]
\centering
\hspace{-.2cm}
\includegraphics[width=0.49\textwidth]{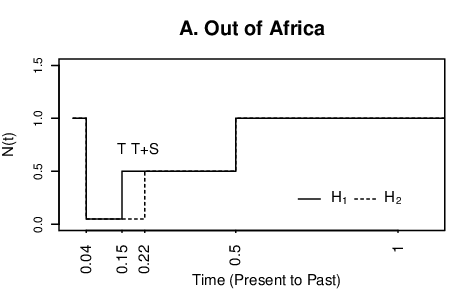}
\vspace{-0.5cm}
\includegraphics[width=0.51\textwidth]{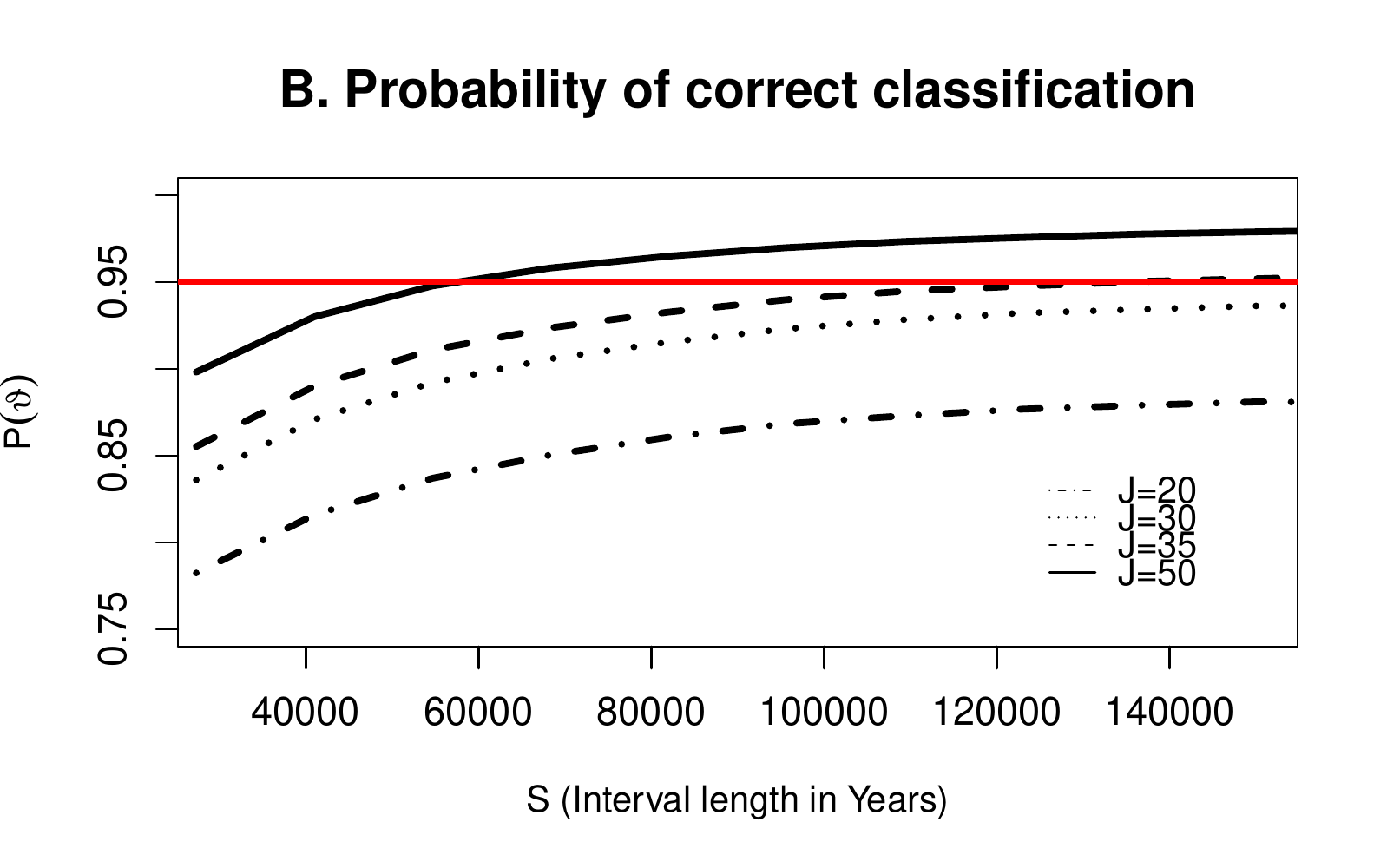}
 \caption{ \textbf{A.} Human population history in coalescent units compatible with previous findings from whole genomes \citep{Li:2011ez}. \textbf{B.} Probability of correct classification $\P[\vartheta(X^J) = \vartheta]$ as a function of the interval length $S$ in years for several values of $J$ loci corresponding to the two hypotheses depicted in A. Red line indicates probability of 0.95.}
\label{fig:MinS}
\end{figure}

Figure \ref{fig:MinS}A shows a population history compatible with a human population history recovered from autosomal DNA in standard coalescent units \citep{Li:2011ez}. In order to convert coalescent parameters into real time and size, time and $N(t)$ need to be divided by the mutation rate per generation. Times need be further multiplied by the generation time.
To make our results comparable to previous studies \citep{Li:2011ez,kim2015can}, we will assume a generation time of $25$ years and that effective population size is expressed in units of $2.732\times 10^{4}$. That is, one unit in the y-axis of Figure \ref{fig:MinS}A corresponds to $2.732 \times 10^{4}$ and one unit in the x-axis of the same plot corresponds to $ 68.3\times 10^{4}$ years. In our analysis, we compare a population trajectory whose second bottleneck starts at time $T=102.45$kya  ($0.15$ in standard units) versus a population trajectory whose second bottleneck starts earlier at time $T+S$ with $S$ ranging from $30$kyr to $150$kyr. Our results from theorem \ref{thm:BayesErrorJLoci} are depicted in Figure \ref{fig:MinS}B. In order to correctly differentiate between the two hypotheses with $S=130$kyr with probability of at least 0.95, we need at least $35$ loci. A correct classification with probability of at least 0.95 is achievable with at least $50$ loci when $S=60$kyr, that is, when the bottleneck started around 162kya versus 102kya.

Our results differ from previously published bounds based on coalescent Bayes error rates.  
\citet{kim2015can} indicate that the minimal $J$ such that any classifier can distinguish between $H_1$ and $H_2$ with probability at least 0.95 and $S=130$kyr is at least $J=10$; while for $S=60$kyr it is $J\approx 20$. A detailed analysis of the differences between our results and previously published bounds of \citep{kim2015can} -- which reflect the fact that we give exact expressions instead of upper bounds on $\P[\vartheta(X^J) = \vartheta]$ -- can be found in section \ref{comparisonKim}.

\subsection{Value of incorporating ancient samples}

\begin{figure}[h]
\centering
 \includegraphics[width=0.9\textwidth]{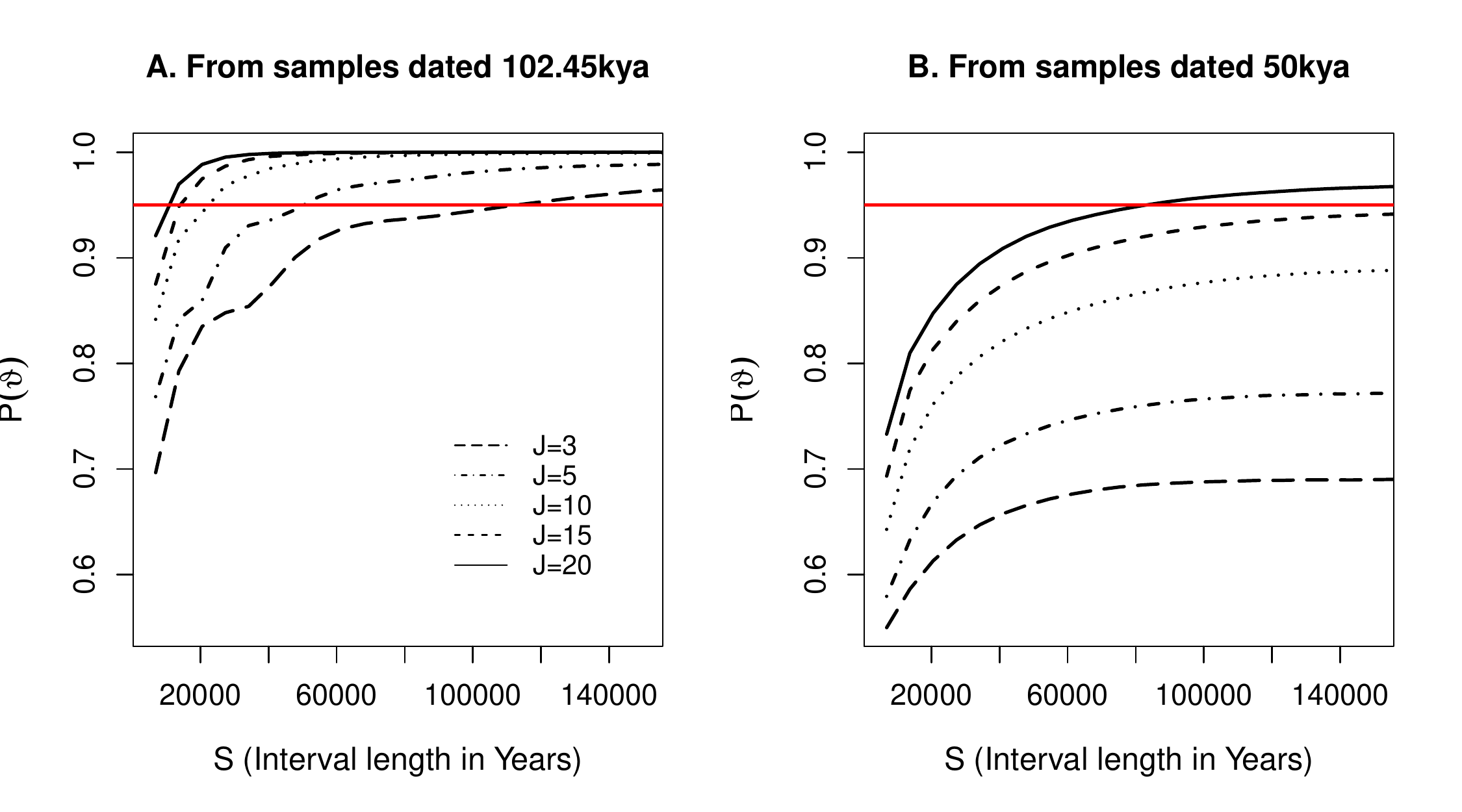}
 \caption{\textbf{Value of incorporating ancient samples.} Probability of correct classification $\P[\vartheta(X^J) = \vartheta]$ as a function of $S$ in years for the human bottleneck example when samples are obtained before the bottleneck (A) and when samples are obtained around 50kya (B) .Curves for different number of loci (J) is indicated by the line patterns. Red line indicates probability of 0.95.
} \label{fig:MinST0}
\end{figure}

Thus far, we have not considered the effect of incorporating samples at different sampling times, and have implicitly assumed that all samples are obtained at present. Results from Theorems \ref{thm:BayesError1Locus} and \ref{thm:BayesErrorJLoci} can be directly applied for the case when the two samples are obtained some time in the past.  In particular, we assess the change in $\P[\vartheta(X^J) = \vartheta]$ when samples are obtained at the end of the bottleneck event at 102.45kya and when samples are obtained at 50kya. These scenarios are  equivalent to putting $\Lambda(T) = 0$ and $\Lambda(T)=1.54$, respectively. Ancient DNA (aDNA) corresponding to $T=50$kya is available from ancient genomes \citep{EuroAge}. 
Obtaining data from the population immediately after the end of the event of interest is in some sense the optimal strategy for statistical inference on that event, and can have an enormous positive effect on inference. This is made clear by Figure \ref{fig:MinST0}A, which shows $\P[\vartheta(X^J) = \vartheta]$ for $J = 2,3,5,10,15$. For all but $J=2$, $\P[\vartheta(X^J) = \vartheta] \ge 0.95$ can be achieved for $S$ greater than about $115$kyr. For $J=15$, it is possible to achieve $\P[\vartheta(X^J) = \vartheta] \ge 0.95$ with $S$ larger than about $15$kyr. When the samples are available from 50kya, it is possible to achieve $\P[\vartheta(X^J) = \vartheta] \ge 0.95$ with at least $J=20$ loci.

\section{Increasing the number of samples}
\begin{figure}[h]
 \centering
 \begin{tabular}{cc}
  \includegraphics[height=0.2\textheight]{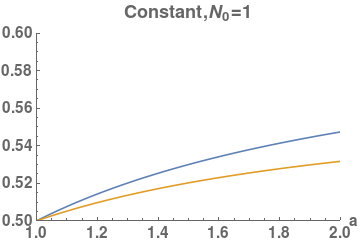} & \includegraphics[height=0.2\textheight]{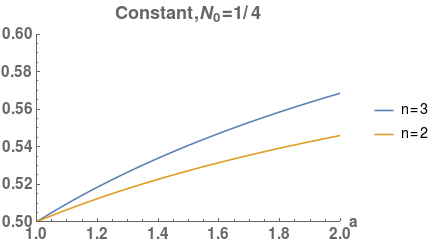} \\
  \includegraphics[height=0.2\textheight]{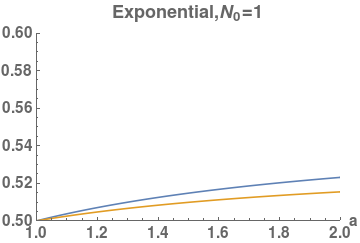} & \hspace{-1.25cm}\includegraphics[height=0.2\textheight]{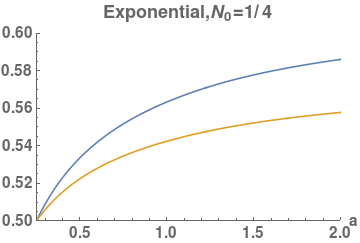} \\
 \end{tabular}
 \caption{\textbf{Effect of adding an additional sample.} Probability of correct classification $\P[\vartheta(X) = \vartheta]$ as a function of $a$ for classification problem \ref{eq:hypo1}. In each case we put $T=1,S=1/2$, and $b=1$. We compare the effect of adding one more sample  ($n=2$  vs $n=3$) for constant and exponential growth population trajectories.
  } \label{fig:Compare23}
\end{figure}

Now we consider the case where $n > 2$. The following theorem gives an exact expression for the probability of correct classification in \eqref{eq:hypo1} from a single locus ($J=1$) when $n=3$.
%
\begin{theorem} \label{thm:BayesTwoTime} Consider the simple hypothesis testing problem of the form \ref{eq:hypo1} when a single genealogy of $n=3$ individuals is observed. Define $\delta$ as in theorem \ref{thm:BayesError1Locus}, then the success rate of the optimal classifier is
 \be
 \P[\vartheta(X) = \vartheta] &= \frac12 + \frac{1}{4} e^{-3 \Lambda(T)} \xi(a,b,N_0,T,S) \\
 \ee
 where
 \be
 \xi(a,b,N_0,T,S) &= 3 e^{2 \Lambda(T)} \left[ e^{-\frac{\delta \wedge S}{a N_0}}-e^{-\frac{\delta \wedge S}{b N_0}}\right] \\
 &+3 \left[ e^{-\frac{2 \left(0 \vee \left(\frac{\delta-S}{2} \wedge S\right)\right)+S}{a N_0}}-e^{-\frac{2  \left(0 \vee \left(\frac{\delta-S}{2} \wedge S\right)\right)+S}{b N_0}} \right]-3 \left[ e^{-\frac{\delta \wedge S}{a N_0}} - e^{-\frac{\delta \wedge S}{b N_0}} \right] \\ 
 &+ \begin{cases} 0 & S < \frac23 \delta \\ 
e^{-\frac{2 \delta}{b N_0}} \left(1+\frac{2 \delta-3S}{b N_0} \right)  - e^{-\frac{2 \delta}{a N_0}} \left(1+\frac{2 \delta-3S}{a N_0} \right) & \frac23 \delta < S < 2\delta \\
e^{-\frac{2\delta}{a N_0}}\left( \frac{4 \delta}{a N_0} + 2 \right) + e^{-\frac{2\delta}{bN_0}}-3e^{\frac{T-2\delta -S}{b N_0}} + e^{-\frac{2 \delta + S}{b N_0}} \frac{3T-3S-4\delta }{b N_0} & S>2 \delta \\
    \end{cases}
 \ee
\end{theorem}
The proof is located in Appendix \ref{sec:BayesTwoTimeProof}. We can use this result to assess how much an additional sample helps in identifying the true population size history. Figure \ref{fig:Compare23} shows four examples of $\P[\vartheta(X) = \vartheta]$ as a function of $a$ while fixing $b=1$; increasing $a$ is equivalent to increasing the separation between the two hypothetical population size histories. In two of the examples, $\Nt(t) = 1$ outside the interval $[T,T+S]$, and in the other two $\Nt(t) = e^t$ outside this interval. In both cases, the probability of identifying the true effective population size function is considerably higher with $n=3$ than $n=2$ when $|a-b|$ is not too close to zero. Thus, additional coalescent times can help considerably to distinguish between alternative histories.

\begin{figure}[h]
 \centering
 \includegraphics[width=0.8\textwidth]{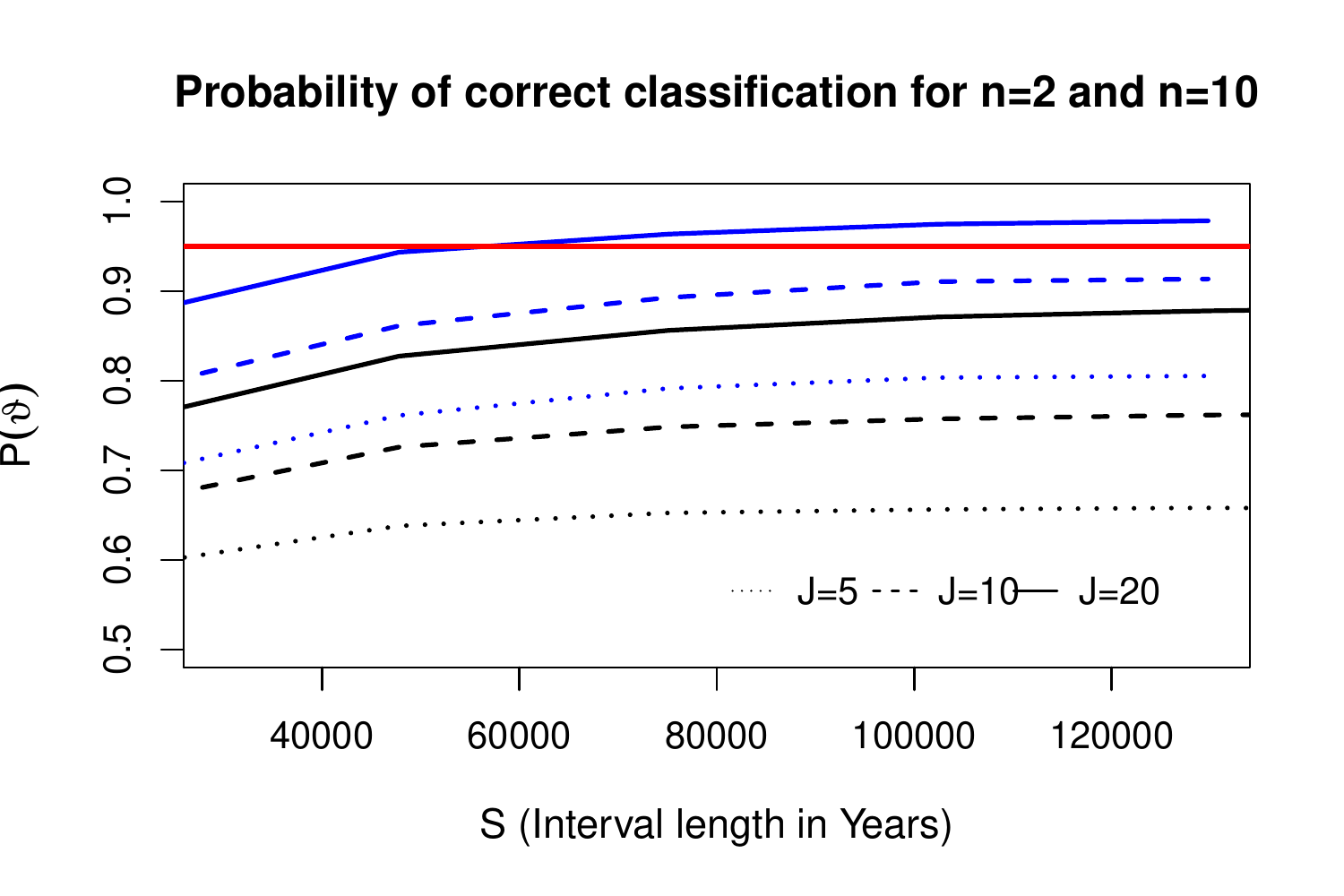}
 \caption{\textbf{Effect of adding more samples in the out-of-Africa scenario} Blue lines represent $\P[\vartheta(X) = \vartheta]$ as a function of  bottleneck length $S$ for $n=10$ and black lines represent  $\P[\vartheta(X) = \vartheta]$ for $n=2$. Different number of loci are distinguished by line patterns. Red line indicates probability of 0.95.} \label{fig:OoAnJ}
\end{figure}

It is clear from the proof of Theorem \ref{thm:BayesTwoTime} that while it is possible to obtain exact expressions for $n>3$, the number of cases that must be treated will grow exponentially in $n$. Of course, it is still possible to approximate $\P[\vartheta(X) = \vartheta]$ by simulation for arbitrary $n$. Here, we re-analyze the out-of-Africa classification problem considered in Section \ref{sec:OoA1} for $n=10$ and $J=1,5,10,20$ as a function of interval length $S$ and compare to $n=2$. The value of $\P[\vartheta(X^J) = \vartheta]$ is approximated by 10,000 Monte Carlo samples. Results are shown in Figure \ref{fig:OoAnJ}. In contrast to the case of $n=2$, where $J=50$ was required to achieve $\P[\vartheta(X^J) = \vartheta]=0.95$ for  $S=60$Kyr, when $n=10$ it is possible to achieve the same success probability with $J=20$. 
Thus, increasing the number of contemporaneous sequences or loci gives sharper inference on the duration of the expansion out-of-Africa.

\section{Bayes error rates in the sequentially Markov coalescent}
We now consider the same classification problem as in \eqref{eq:hypo1} from $J$ consecutive loci along a chromosomal region. We assume the ideal scenario in which we observe the $J$ pairwise coalescent times ($n=2$) corresponding to each of these $J$ loci separated by $J-1$ recombination events. Further, we assume that effective population size trajectories under $H_{1}$ and $H_{2}$ are piece-wise constant functions over time such as the out-of-Africa scenario in Figure \ref{fig:MinS}A. We then approximate 
$\P[\vartheta(X^J)=\vartheta]$ by Monte Carlo from 10,000 simulations generated from each hypothesis under the two coalescent models with recombination: SMC \eqref{eq:smc} and SMC' \eqref{eq:smcprime}. 

We re-analyze the out-of-Africa classification problem considered in section \ref{sec:OoA1} for $n=2$ and $J=2,5,10,20,30,35$ as function of interval length $S$ under independent loci (\ref{thm:BayesErrorJLoci}), SMC' (\ref{eq:smcprime}) and SMC (\ref{eq:smc}). We show that either under SMC' or independent loci, $\P[\vartheta(X^J) = \vartheta]=0.95$ is achievable with $n=35$ loci. For $J<20$, the Bayes error rate in SMC is smaller than the other two alternatives. The significance of this is that it is not necessary to have many independently segregating loci to make inference on features of the historical population size. Instead, virtually the same number of non-independent loci separated by recombination events will suffice. The set of all dependent loci is of course considerably larger than the largest set of independent loci, so the result suggests optimism in the potential to reconstruct features of the population size trajectory in the relatively distant past.

\begin{figure}[h]
 \centering
 \includegraphics[width=1.0\textwidth]{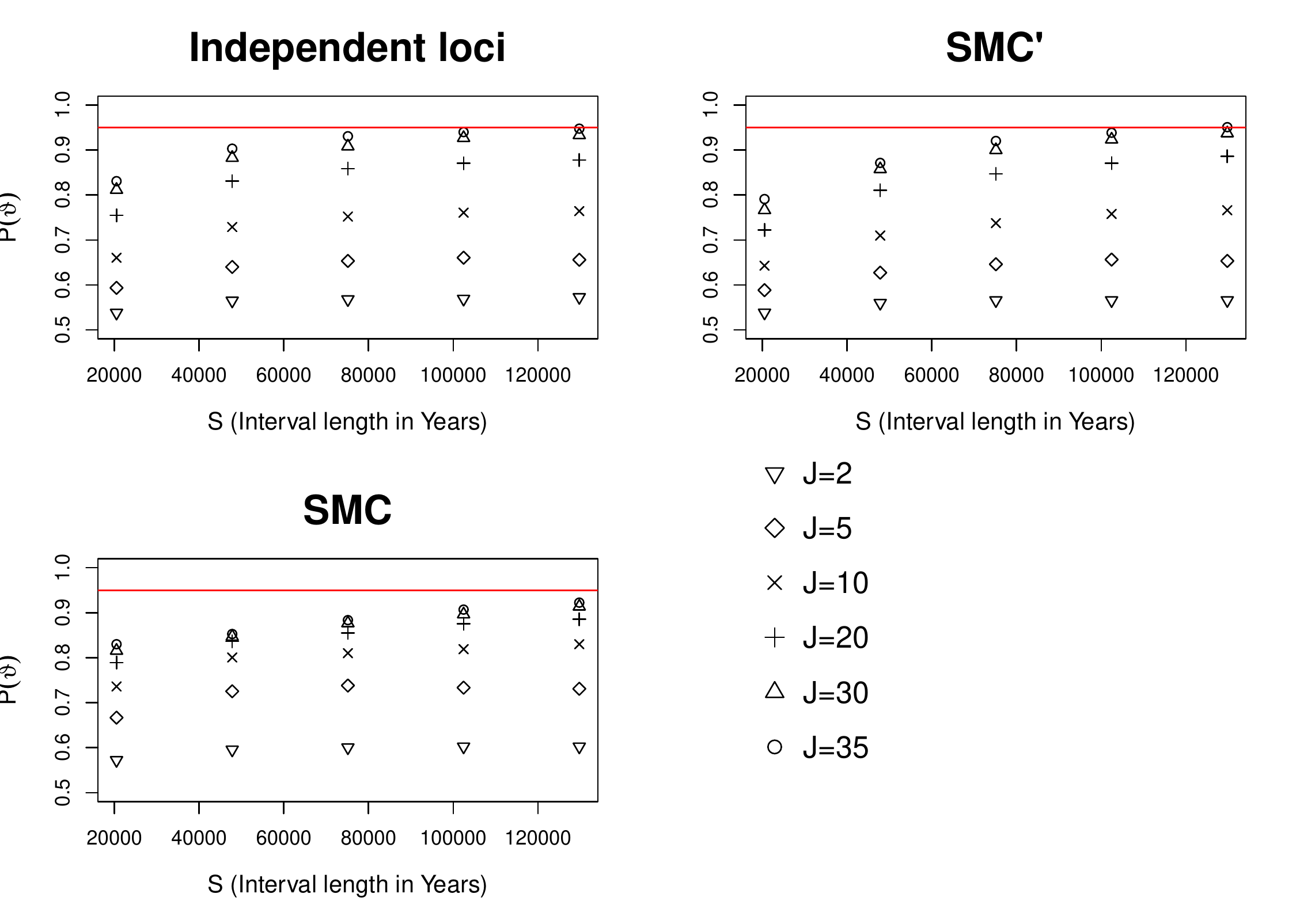}
 \caption{\textbf{Sequentially Markov coalescent in the out-of-Africa scenario} Probability of correct classification under independent sampling, SMC' and SMC. Different patterns represent different number of loci. Red line indicates probability of 0.95.} \label{fig:smcres}
\end{figure}

\section{Other scenarios}
We now consider a more general classification problem when pairwise coalescent data is available at a single locus or multiple loci :
\be \label{eq:hyp2}
H_1 : \Nt(t) = N_1(t) \\
H_2 : \Nt(t) = N_2(t).
\ee
We consider the case where $N_2(t) = c N_1(t)$ for $c \in (0,1)$, where analytic expressions for $\P[\vartheta(X) = \vartheta]$ are available even when $N_1(t)$ is not piecewise constant.
\begin{theorem} \label{thm:BayesError2_1Locus}
Consider the simple hypothesis testing problem of the form \eqref{eq:hyp2} such that $N_2(t) = c N_1(t)$ with $0 < c < 1$ when a single pairwise coalescent time is observed ($n=2$) and assign equal prior probabilities to both hypotheses. Then the probability of correct classification is:
 \be
 \P[ \vartheta(X) = \vartheta] &= \frac12 c^{\frac{c}{1-c}} + \frac12 \left(1 - c^{\frac{1}{1-c}} \right).
 \ee
\end{theorem}

%

\begin{theorem} \label{thm:BayesError2_JLoci}
Consider the conditions of Theorem \ref{thm:BayesError2_1Locus} for $J$ independent loci. The probability of correct classification is
\be \label{eq:ConsShiftJ}
\P[\vartheta(X^J) = \vartheta] &= \frac12 \left( 1 - \frac1{\Gamma(J)} \left[ \gamma\left( J,\frac{-J c \log c}{1-c} \right) - \gamma\left( J,\frac{J \log c}{c-1} \right) \right] \right). 
\ee
where $\gamma(a,b)$ is the lower incomplete gamma function.
\end{theorem}

\begin{figure}[h]
\centering
 \includegraphics[width=0.4\textwidth]{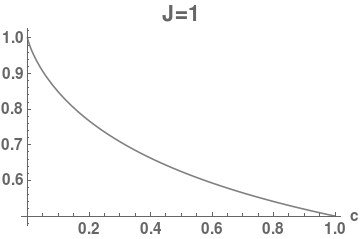}
\includegraphics[width=0.5\textwidth]{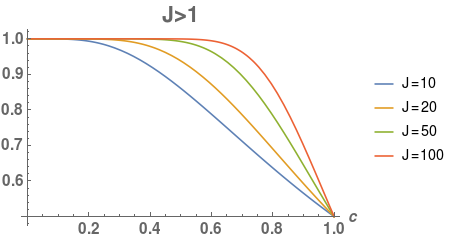}
\caption{\textbf{Left.} Probability of correct classification $\P[\vartheta(X) = \vartheta]$ as in \eqref{eq:ConsShiftJ} when $N_2(t) = c N_1(t)$ and $J=1$ (Theorem \ref{thm:BayesError2_1Locus}). \textbf{Right. }$\P[\vartheta(X^J) = \vartheta]$ as a function of $c$ for several values of $J$. } \label{fig:ConsShiftProb}
\end{figure}

Figure \ref{fig:ConsShiftProb} shows \eqref{eq:ConsShiftJ} as a function of $c$ for different values of $J$. As expected, the larger $J$, the larger the value of $c$ at which high probability of identifying the true population size history can be achieved. However, even for $J=100$, we must have $c \approx 0.75$ or smaller to give probability 0.95 of selecting the true population size history.
Our results from Theorems \ref{thm:BayesError2_1Locus} and \ref{thm:BayesError2_JLoci} differ from previously published Bayes error rates bounds \citep{kim2015can}. In the following section, we present a more detailed analysis of the differences between our exact expressions and the bounds \citep{kim2015can}.

\section{Bounds on Bayes classification error.}  \label{comparisonKim}
\citet{kim2015can} provided lower bounds on  Bayes error rates from pairwise coalescent data. We now provide a comparison of some of our results to these previously published bounds. In this section, we will let $Y$ denote a random coalescence time generated under $H_1$ and $Z$ denote a random coalescence time generated under $H_2$. Assuming a classification problem of the form \eqref{eq:hypo1} and prior probability 1/2 on $H_1$ 
and $H_2$, 
the Bayes error rate for any classifier is at least $(1-\Upsilon)/2$ where
\be
\Upsilon = 
\d_{\TV}(Y,Z)=\d_{\TV}(\mu, \nu) \equiv \sup_A | \mu(A) - \nu(A) | 
\ee
is the total variation distance between probability measures $\mu,\nu$, such that $Y \sim \mu, Z \sim \nu$. 
%
The authors then apply the inequality 
\be
\frac12 \d_{\TV}^2 \le \d_{\H}^2
\ee
where $\d_{\H}$ is the \emph{Hellinger distance}. Let $P$ and $Q$ be 
probability measures that are absolutely continuous with respect to some 
dominating measure $\lambda$, and let $f_P = \frac{dP}{d\lambda}$, 
$f_Q=\frac{dQ}{d\lambda}$ be their respective Radon-Nikodym derivatives. The 
Hellinger distance between $P$ and $Q$ is defined by
\be
\d^2_{\H}(P,Q) = \frac12 \int (\sqrt{f_P} - \sqrt{f_Q})^2 d\lambda.
\ee
In the case where $\lambda$ is Lebesgue measure, $f_P$ and $f_Q$ are the 
densities of $P$ and $Q$. The main result of \cite{kim2015can} is
\begin{theorem}[\citet{kim2015can}, Theorem 1] \label{thm:Hellinger}
Suppose $n=2$ in \eqref{eq:IsoLikelihood}. Then
\be
\d_{\H}^2(Y,Z) = e^{-\int_0^T \frac{1}{\Nt(t)}dt} 
\left(1-e^{-\frac{(a+b) S}{2abN_0}} \right) \frac{(\sqrt{a} - 
\sqrt{b})^2}{a+b}.
\ee
\end{theorem}
We give a proof in the appendix that fills in some additional details of the 
proof in \cite{kim2015can}. Rather than obtaining bounds on the Bayes error rate using the Hellinger distance,
we compute the probability of correct inference on $\vartheta$.

\begin{figure}[h]
\centering
\begin{tabular}{cc}
\includegraphics[height=0.2\textheight]{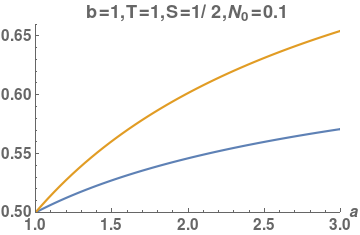} & \includegraphics[height=0.2\textheight]{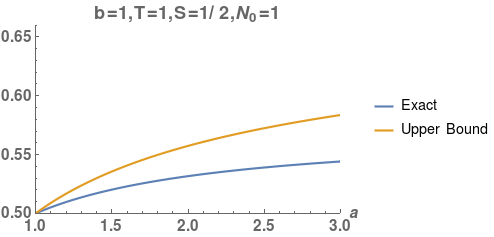} \\
\includegraphics[height=0.2\textheight]{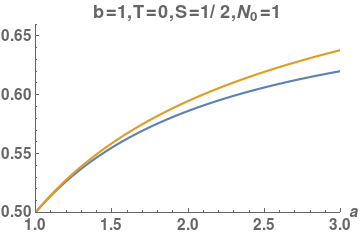} & \hspace{-2.45cm} \includegraphics[height=0.2\textheight]{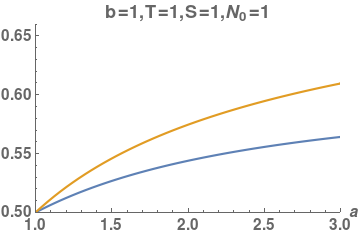} \\
\end{tabular}
\caption{Exact $\P[\vartheta(X) = \vartheta]$ (blue) and upper bound on this quantity from \cite{kim2015can} (yellow) for different values of $T,S,$ and $N_0$.   } \label{fig:CompareBounds}
\end{figure}

In Figure \ref{fig:CompareBounds}, we compare our results to the Hellinger bounds of \cite{kim2015can} for different values of $a,b, N_0$. 
The upper bound based on the Hellinger distance from \cite{kim2015can} is given by 
\be
\frac12 + \frac12 \sqrt{2 H^2(f_1,f_2)}
\ee
with $H^2(f_1,f_2)$ as in \eqref{eq:HellingerBoundOneTime}. Evidently the
Hellinger bound is quite loose when $|a-b|$ is not near zero.

\citet{kim2015can} use the inequality 
\be \label{eq:HellingerIneq}
\d^2_H(Y^J,Z^J) \le J \d^2_H(Y,Z),
\ee
which holds when the $J$ genealogies are independent, in combination with Theorem \ref{thm:Hellinger}, to obtain lower bounds on the error rate for $J$ independent loci. They then use these lower bounds to calculate quantities like bounds on the minimal $S$ such that the correct hypothesis will be selected with probability 0.95 for several examples. 
In contrast, our results give the exact value of $\P[\vartheta(X^J) = \vartheta]$,
which allows us to compute exactly the value of $S$ to achieve the desired  Bayes error rate for any $J$. The results on the minimal number of loci $J$ necessary to achieve a fixed error rate differ substantially from the results in \cite{kim2015can}. The looseness of the bound on $\P[\vartheta(X^J) = \vartheta]$ obtained using the Hellinger distance 
is clear from Figure \ref{fig:CompareBounds}.

%

\begin{figure}[h]
 \centering
 \begin{tabular}{cc}
  \includegraphics[height=0.2\textheight]{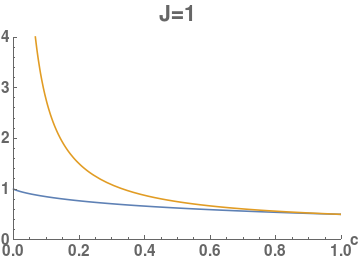} & \includegraphics[height=0.2\textheight]{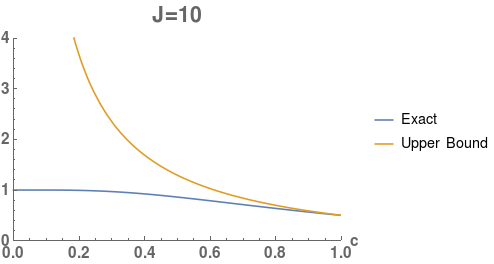}
 \end{tabular}
\caption{ Exact $\P[\vartheta(X^J) = \vartheta]$ compared to the upper bound from Theorem 3.2 of \cite{kim2015can} as a function of $c$ for two different values of $J$.  } \label{fig:ConsShiftCompare}
\end{figure}

Moreover, the expression in \eqref{eq:ConsShiftJ} can be directly compared with Theorem 3.2 of \cite{kim2015can}. Translated into our notation and conventions, this result states that
\be \label{eq:KimConsShift}
\P[\vartheta(X^J) = \vartheta] \le \frac12 + \frac14 \sqrt{J(n-1)} \left( \frac1c - 1 \right). 
\ee
Figure \ref{fig:ConsShiftCompare} shows the bound from \eqref{eq:KimConsShift} along with the exact probability of identifying the true $\Nt(t)$ as a function of $c$ for $n=2$ and different values of $J$. The bound is apparently quite loose when $c$ is not close to 1. It becomes trivial (greater than 1) for $c \approx 0.4$ when $J=1$ and $c \approx 0.7$ when $J=10$. 

\section{Bayes risk under conjugate priors}
\begin{figure}[h]
 \centering
 \includegraphics[width=0.5\textwidth]{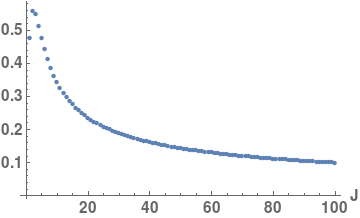}
 \caption{Root risk of Bayes estimator with $\alpha=\beta=1$ and $c=1$. } \label{fig:RiskBayes}
\end{figure}

Although our focus has been on inferential limits for distinguishing among two states of nature, we briefly consider estimation of a constant population size trajectory. We asses the risk of estimators of the function $\Lambda(x)$ in the case of $n=2$ and $\Nt(t) = \frac1c$ with conjugate priors on $c$. In this setting, the coalescent time $x$ is $\text{Exponential}(c)$ with conjugate prior $c \sim \Gam(\alpha,\beta)$ and for a sample of $J$ independent pairwise coalescent times we have
\be
c \mid x^{(1)},\ldots,x^{(J)} \sim \Gam(\alpha+J,\beta + J \bar x)
\ee
with posterior expectation
\be
\E[c \mid x^{(1)},\ldots,x^{(J)}] &= \frac{\alpha+J}{\beta+J \bar x}.
\ee
Note that
\be
J \bar X \mid c \sim \Gam(J,c)
\ee
so the squared error risk of the Bayes estimator of $c$ is
\be
R(\hat c, c) := \int_0^{\infty}  \left( \frac{\alpha+J}{\beta+z} - c\right)^2 \frac{\beta^\alpha}{\Gamma(\alpha)} z^{J-1} e^{-c z} dz.
\ee
This can be evaluated in terms of analytic functions for any $\alpha,\beta$, but for simplicity, we assume $\alpha = \beta = 1$, so that
\be
R(\hat c, c) &= \int_0^{\infty}  \left( \frac{1+J}{1+z} - c\right)^2 e^{-z} dz \\
&= c \left(-(J+1) e^c \left(J^2+(J+3) c-1\right) E_J(c)+(J+1)^2+c\right)
\ee
where 
\be
E_J(c) = \int_1^{\infty} \frac{e^{-cz}}{z^J} dz, 
\ee
is a generalization of the exponential integral function. Figure \ref{fig:RiskBayes} shows the square root of risk as a function of the number of loci $J$ for values of $J \in \{1,\ldots,100\}$ with $c=1$. The root risk  decreases logarithmically in $J$; it is approximately 0.1 for $J=100$, and about 0.24 for $J=20$. Thus, if one wants the root risk to be small relative to the truth, it is necessary to have $J$ rather large. In this example, in order to have the root risk be about 10 percent the magnitude of the truth, we need $J \approx 100$.

\section{Discussion}


Availability of ancient and present-day DNA samples from a population allows statistical reconstruction of the effective population size trajectory. The effective population size is a measure of relative genetic diversity whose actual magnitude is not easily interpreted in units of census population size \citep{Wakeley10112008}. However, changes of effective population size over time are informative about the genetic history of the population. In this manuscript, we assess the ability to differentiate or classify between alternative hypotheses about the effective population size. 

Assessment of inferential limits 
in population genetic studies is becoming important in the face of ongoing large-scale studies of genetic variation. Statistical methods are usually restricted to small samples or rely on approximations and insufficient summary statistics. As such, choosing the optimal subset of data with which to perform statistical inference is of great interest. Aspects of the data and adequacy of the model will affect the ability to draw meaningful conclusions. Here, we have eliminated the effect of factors such as data quality, sample selection and sequence alignment and concentrated on the ideal scenario of having a complete realization of the genealogical process free of errors. In practice, genealogies are not available and instead we observe DNA sequence variation; therefore our results are upper bounds on the achievable probability of recovering the true population size history in population genetic studies. These results provide guidance to practitioners in choosing a sampling design subject to computational constraints. In particular, they give insight into the key questions of which scientific hypotheses can be assessed, and the optimal choice of the number of loci, sampling times, and individuals to include in a sample to achieve a specific inferential goal. They also offer a possible explanation for disagreement in the literature over timing and duration of historical genetic events such as the out-of-Africa human population bottleneck, suggesting that some studies may simply not have sufficient data to distinguish between the hypotheses of interest with high probability.

\citet{Fu693,PluzhnikovDonnelly}, and \citet{Felsenstein2006} argued that in the constant population model ($\theta=2N\mu$), accuracy of estimators of $\theta$ 
increases linearly in the number of independent loci, logarithmically in the number of samples, and is unaffected by sequence length. In the coalescent with variable population size, \citet{MYERS2008342} showed that estimators based on the SFS cannot distinguish between two alternative hypotheses. \citet{Terhorst23062015} showed that estimators of $N(t)$, based on the same statistic SFS, have minimax rate of convergence that is logarithmic in the number of independent loci and independent of the number of individuals sampled. Our results support a more complex view of the value of additional samples or loci. While in general, the improvement in the probability of recovering the true population history appears to be sublinear in both $J$ and $n$, the improvement from adding an additional sample or locus depends greatly on the details of the two hypotheses being considered. For example, increasing from $n=2$ to $n=3$ samples can in some cases double the excess probability of recovering the truth $\P[\vartheta(X) = \vartheta]-1/2$ (the probability is always lower bounded by 1/2). In general, smaller improvements are seen from increasing $J$, but we have demonstrated that high probability of recovering the true population size history in the out-of-Africa example is attainable using values of $n$ and $J$ that are available from modern datasets and for which exact computation is feasible. In addition, our  
results suggest that incorporation of ancient genomes is the optimal strategy to improve inferential performance  
in the expansion out-of-Africa problem, which is of significant interest in human population genetics.

\citet{PluzhnikovDonnelly} considered the constant population model with recombination and argued that when the recombination rate is high, increasing the sequence length effectively increases the number of independent loci. Indeed, when two genomic segments are separated by a recombination event, individuals at these two segments (loci) derive from two different but correlated genealogies. As the number of recombination events increases, the correlation between the two genealogies becomes weaker, and hence increasing the length of sequenced segments increases the opportunity to observe a larger number of realizations from genealogically independent loci \citep{Palaciosgenetics,griffiths1997ancestral}. Our results for the pairwise SMC' model of recombination support the remarkable conclusion that loci separated by recombination events have nearly the statistical value as the same number of independent loci. This suggests that pairwise SMC' is a very powerful framework for inference of population size trajectories. An interesting future direction is to explore the effect of increasing the number of individuals in the SMC'.

\newpage
\appendix

\section{Proof of theorem \ref{thm:BayesError1Locus}}
\begin{proof}
Define
\be
\Lambda(w,x) \equiv \int_w^x \frac{1}{\Nt(t)} dt.
\ee
For shorthand we write $\Lambda(x) = \Lambda(0,x)$. 
$\Lambda : \mathbb R_+ \to \mathbb R_+$ is a monotone strictly increasing 
function, which is enough to guarantee the existence of an inverse
\be
\Lambda^{-1}(t) = x \Longleftrightarrow \Lambda(x) = t,
\ee
The likelihood ratio for $H_1$ vs $H_2$  (\ref{eq:hypo1}) can be expressed by
\be \label{eq:DensityOnetime}
\log \BF_{12}(x) = \left\{ \begin{array}{lc} 0 & 
x<T \\ 
\log \frac{b}{a} -\frac{x-T}{a N_0} + \frac{x-T}{b N_0} & T \le x < T+S \\
 \frac{S}{b N_0}-\frac{S}{a N_0}  & T+S \le x 
\end{array} \right.
\ee
Notice that the waiting time until the coalescent event has survival function
\be
\mathbf P[ X > x] = e^{-\Lambda(x)}.
\ee
Now we want to calculate $\mathbf P[\vartheta(X) = 1 \mid H_1]$. Assume that 
if $\log \BF_{12}(x) = 0$ we select either $H_1$ or $H_2$ by flipping a fair 
coin. If $a>b$ then 
\be
\log \BF_{12}(x) > 0, T \le x \le T+S \Longleftrightarrow x > \delta+T, 
\ee
and if $b > a$
\be
\log \BF_{12}(x) > 0, T \le x \le T+S \Longleftrightarrow x < \delta+T .
\ee

Assuming $a>b$ and denoting $f_{i}(x)$ the density under $H_{i}$ for $i=1,2$, we have 
\be
\mathbf P[\vartheta(X) = 1 \mid H_1] &= \frac12 \mathbf P[X < T \mid H_1] + 
\int_T^{T+S} \bone\left\{x > \frac{ab}{b-a} N_0 \log \frac{b}{a} + T \right\} 
f_1(x) dx \\
&+ \bone\{ b < a\} \mathbf P[X > T+S \mid H_1] \\
&= \frac12 (1-e^{-\Lambda(T)})  + 
\int_{T+(\delta \wedge S)}^{T+S} e^{-\Lambda(T)} \frac{1}{aN_0} e^{-\frac{x-T}{a N_0}} dx + \bone_{\{b < a\}} e^{-\Lambda(T)-\frac{S}{a N_0}}  \\
&= \frac12 (1-e^{-\Lambda(T)}) + e^{-\Lambda(T)} \left[e^{-  \frac{\delta \wedge S}{a N_0}} - e^{-\frac{S}{a 
N_0}} \right]  + \bone_{\{b < a\}} e^{-\Lambda(T)-\frac{S}{a N_0}} ,
\ee  
and
\be
\mathbf P[\vartheta(X) = 2 \mid H_2] &= \frac12 \mathbf P[X < T \mid H_2] + 
\int_T^{T+S} \bone\left\{x < \frac{ab}{b-a} N_0 \log \frac{b}{a} + T \right\} 
f_2(x) dx  \\
&+ \bone\{ a < b\} \mathbf P[X > T+S \mid H_2] \\
&= \frac12 (1-e^{-\Lambda(T)})  + 
\int_T^{T+(\delta \wedge S)} e^{-\Lambda(T)} \frac1{b N_0} e^{-\frac{x-T}{b N_0}} dx  + \bone_{\{a < b\}} e^{-\Lambda(T)-\frac{S}{b N_0}}  \\
&= \frac12 (1-e^{-\Lambda(T)}) + e^{-\Lambda(T)} \left[1 - e^{-  \frac{  
\delta \wedge S}{b N_0}} \right] + \bone_{\{a < b\}}
e^{-\Lambda(T)-\frac{S}{b N_0}}. 
\ee
Assuming equal prior probability of $H_1$ and $H_2$ we get
\be
\mathbf P[\vartheta(X)=\vartheta] &= \frac12 (1-e^{-\Lambda(T)} 
+ e^{-\Lambda(T)-\frac{S}{(a \vee b) N_0}}) \\
&+ \frac12 e^{-\Lambda(T)} \left[ e^{-\frac{ \delta \wedge S}{aN_0}}- e^{- \frac{S}{a N_0}} 
 \right]  + \frac12 e^{-\Lambda(T)} \left[ 1 - e^{- \frac{\delta \wedge S}{b N_0}} \right] \\
 &= \frac12 + \frac12 e^{-\Lambda(T)} \left( e^{-\frac{S}{(a \vee b) N_0}}-e^{-\frac{S}{a N_0}} \right) + \frac12 e^{-\Lambda(T)} \left( e^{-\frac{ \delta \wedge S}{aN_0}} - e^{- \frac{\delta \wedge S}{b N_0}} \right) \\
  &= \frac12 + \frac12 e^{-\Lambda(T)} \left( e^{-\frac{ \delta \wedge S}{aN_0}} - e^{- \frac{\delta \wedge S}{b N_0}} \right).
\ee
This assumed $a>b$. If instead $b>a$ then the inequalities in the integrand 
when we integrate between $T$ and $T+S$ would be reversed, so the exact expression for any $a>0,b>0$ is
\be \label{eq:ExactP}
\mathbf P[\vartheta(X)=\vartheta] 
 &= \frac12 + \frac12 e^{-\Lambda(T)} \left( e^{-\frac{ \delta \wedge S}{(a \vee b)N_0}} - e^{- \frac{\delta \wedge S}{(a \wedge b) N_0}} \right).
\ee
\end{proof}

\section{Proof of Theorem \ref{thm:BayesErrorJLoci}}
\begin{proof}
Fix an integer $J \ge 1$ and define $a_1 = a, a_2 = b$ for ease of notation. 
We first define the following auxiliary functions
\be
Q_{i}(T) &\equiv e^{-\int^{T}_{0}\frac{dt}{N_{i}(t)}}, \quad Q_{i}(T,T+S) 
\equiv e^{-\int^{T+S}_{T}\frac{dt}{N_{i}(t)}} \\
q_{i}(T) &\equiv \frac{1}{N_{i}(T)}e^{-\int^{T}_{0}\frac{dt}{N_{i}(t)}}, \quad
q_{i}(T,T+S) \equiv 
\frac{1}{N_{i}(T+S)}e^{-\int^{T+S}_{T}\frac{dt}{N_{i}(t)}}
\ee
The coalescent density for a coalescent time with effective population size trajectory $N$ for the intervals $(0,T]$ and $(T+S,\infty)$ and $N_{i}$ for the interval $(T,T+S]$ is
\be
f_{i}(t)=\begin{cases}
q(t) & 0<t<T\\
Q(T)q_{i}(T,t) & T\leq t<T+S\\
Q(T)Q_{i}(T,T+S)q(T+S,t) & t \geq T+S\\
\end{cases}
\ee
so that the likelihood ratio for a single time point can be 
expressed as
\be
\frac{f_1(x^j)}{f_2(x^j)}& =\left[ 
\frac{q_1(T,x^j)}{q_2(T,x^j)}\right]^{\bone\{T \leq x^j < T+S\}} 
\left[\frac{Q_{1}(T,T+S)}{Q_{2}(T,T+S)}\right]^{\bone\{x^j \geq T+S\}}\\
& = \left[\frac{b}a e^{-\frac{(b-a)(x^j-T)}{ab N_0}}\right]^{\bone\{T \leq  
x^j < T+S\}}\left[e^{-S\frac{(b-a)}{ab N_0}}\right]^{\bone\{ x^j \geq T+S \}},
\ee
giving 
\be
\log \prod^J_{j=1} \frac{f_1(x^j)}{f_2(x^j)} &=  
\sum^J_{j=1}\bone\{T<x^j \leq T+S\} \left[ \log \frac{b}a 
-(x^j-T)\frac{(b-a)}{ab N_0}\right] \\
&- \sum^J_{j=1}\bone\{x^j> T+S\} \frac{S(b-a)}{ab N_0}.
\ee

Defining
\be
\ell_1=\sum^J_{j=1}\bone\{x^j \leq T\}, \quad \ell_2 = 
\sum^J_{j=1}\bone\{T<x^j \leq T+S\}, \quad 
\ell_3=\sum^J_{j=1} \bone\{x^j \geq T+S\},
\ee
we have that $\log \BF_{12} > 0$ when
\be
\sum^J_{j=1}\bone\{T<x^j \leq T+S\} \left[ \log \frac{b}a 
-(x^j-T)\frac{(b-a)}{ab N_0}\right] &> \sum^J_{j=1}\bone\{x^j> T+S\} 
\frac{S(b-a)}{ab N_0} \\
\ell_2 \left(\log \frac{b}a + T\frac{(b-a)}{ab N_0} \right) - \ell_3 S 
\frac{(b-a)}{ab N_0} &> \frac{(b-a)}{ab N_0} \sum_{j : x^j \in [T,T+S]} x^j \\
\ell_2 \left( \frac{ab N_0}{a-b} \log \frac{a}b + T \right) 
- \ell_3 S &< \sum_{j : x^j \in [T,T+S]} x^j, \\
\ell_2 \left( \delta + T \right) 
- \ell_3 S &< \sum_{j : x^j \in [T,T+S]} x^j, 
\ee
where the inequality reversed since $(b-a)/(ab N_0)$ is negative.

Now, $\log \BF_{12}=0$ only if $x^j<T$ for all $j=1,\ldots,J$. In this case, 
we flip a fair coin and accept $H_1$ if it shows heads. Moreover,
if $L_2 = 0$ and $L_3 > 0$, then $\log \BF_{12} > 0$. Denote by 
$L=(L_1,L_2,L_3)$ the random vector whose observed entries are 
$\ell=(\ell_1,\ell_2,\ell_3)$. Notice that for a generic coalescent time $X$
\be
L \mid H_i &\sim \Multi(J,\mathbf{p}) \\
p_1 &= \P[X \le T] = (1-e^{-\Lambda(T)}) \\
p_2 &= \P[T < X \le T+S] = (e^{-\Lambda(T)}-e^{-\Lambda(T)-\frac{S}{a_i N_0}}) 
\\
p_3 &= \P[X > T+S] = e^{-\Lambda(T)-\frac{S}{a_i N_0}},
\ee
and we have
\be
\P[ \vartheta(X^J)=1 \mid H_1] &=\frac12 \P(L_1=J \mid H_1) + \P(L_2 = 0,L_3>0 \mid 
H_1) \\
&+ \sum_{(\ell_2,\ell_3) : \ell_2 > 0} \P(L=\ell \mid H_1)\P(\BF_{12}(X^J)>0\mid 
L=\ell, H_1)
\ee
with
\be
\P[\BF_{12}(X^J) > 0 \mid L=\ell,H_1] &= \P\left[ \sum_{j: X^j \in [T,T+S]} X^j > 
\ell_2 (\delta + T) - \ell_3 S \,\bigg| \, L=\ell \right] \\
&= \P\left[ \sum_{j: X^j \in [T,T+S]} X^j > \ell_2 
(\delta + T) - \ell_3 S \,\bigg| \, T < X^j \le T+S \right] \\
&= \P\left[ \sum_{j=1}^{\ell_2} X^j_* > \ell_2 \delta - \ell_3 S 
\, \bigg|\, X_*^j < S \right]
\ee
for $X^j_*$ independent exponential random variables with rate $(a N_0)^{-1}$. 
So letting $W^*(\ell_2) = \sum_{j=1}^{\ell_2} X_*^j$, the relevant probabilities 
involve the CDF of the sum of $\ell_2$ many independent exponentials with rate 
$(a N_0)^{-1}$ truncated to the interval $[0,S]$, and we have
\be
\P[ \vartheta(X^J) = 1 \mid H_1] &=\frac12 \P(L_1=J \mid H_1) + \P(L_2 = 0,L_3>0 \mid 
H_1) \\
&+ \sum_{(\ell_2,\ell_3) : \ell_2 > 0} \P(L=\ell \mid H_1) \P[W^*(\ell_2) > \ell_2 \delta - \ell_3 S \mid H_1].
\ee
It follows then that since $\P(L_1 = J \mid H_1) = \P(L_1 = J \mid H_2)$, the Bayes error rate can be written as
\be \label{eq:BayesErrorJ}
\P[\vartheta(X^J) = \vartheta ] &= \frac12 \P(L_2=0 \mid 
H_1) \\
&+ \frac12 \sum_{(\ell_2,\ell_3) : \ell_2 > 0} \P(L=\ell \mid H_1) \P[W^*(\ell_2) > \ell_2 \delta - \ell_3 S \mid H_1] \\
&+ \frac12 \sum_{(\ell_2,\ell_3) : \ell_2 > 0} \P(L=\ell \mid H_2) \P[W^*(\ell_2) < \ell_2 \delta - \ell_3 S \mid H_2].
\ee
\end{proof}

\section{Proof of Theorem \ref{thm:BayesTwoTime}} \label{sec:BayesTwoTimeProof}
Recall we are studying the case where $H_1 : N = N_1(t)$ and $H_2 : N = N_2(t)$ and
\be
N_1(t) &= \begin{cases} \Nt(t) & 0 \le t \le T \\ a N_0 & T \le t \le T+S \\ \Nt(t) & t > T+S \end{cases} \\
N_2(t) &= \begin{cases} \Nt(t) & 0 \le t \le T \\ b N_0 & T \le t \le T+S \\ \Nt(t) & t > T+S \end{cases} 
\ee
for $\Nt(t)$ any bounded, strictly non-negative function.

\begin{enumerate}
\item \emph{Case 1: $0 < x_2 < x_1 < T$}. In this case the likelihood under either
$H_1$ or $H_2$ is the same
\be
L(x_2,x_1 \mid \Nt(t))  &= \frac3{\Nt(x_2) \Nt(x_1)} e^{-2 \Lambda(x_2)-\Lambda(x_1)} \\
\ee
and so
\be
\log \BF_{12}(x) = 0.
\ee

\item \emph{Case 2: $0 < x_2 < T < x_1 < T+S$}. In this case the likelihood under $H_i$ is
\be
L(x_2,x_1 \mid \Nt(t)) &= \frac{3}{\Nt(x_2)} \frac{1}{a_i N_0} e^{-2\Lambda(x_2)-\Lambda(T)-\frac{x_1-T}{a_i N_0}} 
\ee
so designating $a_1=a, a_2 = b$ as before
\be
\log \BF_{12}(x) &= \log \frac{b}{a} -\frac{x_1-T}{a N_0} + \frac{x_1-T}{b N_0} \\
&= \log \frac{b}{a} + \frac{(a-b)(x_1-T)}{ab N_0}.
\ee

\item \emph{Case 3: $0 < x_2 < T < T+S < x_1$}. In this case the likelihood under $H_i$ is
\be
L(x_2,x_1 \mid \Nt(t)) &= \frac{3}{\Nt(x_2)} \frac{1}{\Nt(x_1)} e^{-2\Lambda(x_2)-\Lambda(T)-\frac{S}{a_i N_0}-\Lambda(T+S,x_1)} \\
\ee
so 
\be
\log \BF_{12}(x) &= \frac{S}{b N_0} - \frac{S}{a N_0} = \frac{(a-b)S}{ab N_0}.
\ee

\item \emph{Case 4: $0 < T < x_2 < x_1  < T+S$}. In this case the likelihood under $H_i$ is
\be
L(x_2,x_1 \mid \Nt(t)) &= \frac{3}{a_i N_0} \frac{1}{a_i N_0} e^{-3\Lambda(T)-\frac{3 (x_2-T)}{a_i N_0}-\frac{x_1-x_2}{a_i N_0}}  \\
&= \frac{3}{a_i^2 N_0^2} e^{-3\Lambda(T)} e^{-\frac{2x_2+x_1-3T}{a_i N_0}}
\ee
so 
\be
\log \BF_{12}(x) &= 2 \log \frac{b}{a} -\frac{2x_2+x_1-3T}{a N_0} + \frac{2x_2+x_1-3T}{b N_0} \\
&= 2 \log \frac{b}{a} + \frac{(a-b)(2x_2+x_1-3T)}{ab N_0}
\ee

\item \emph{Case 5: $0 < T < x_2 < T+S < x_1$}. In this case the likelihood under $H_i$ is
\be
L(x_2,x_1 \mid \Nt(t)) &= \frac{3}{a_i N_0} \frac{1}{\Nt(x_1)} e^{-2\Lambda(T)-\frac{2 (x_2-T)}{a_i N_0}-\Lambda(T)-\frac{S}{a_i N_0}-\Lambda(T+S,x_1)}  \\
&= \frac{3}{a_i N_0} \frac{1}{\Nt(x_1)} e^{-3\Lambda(T)- \Lambda(T+S,x_1)} e^{-\frac{2 (x_2-T)+S}{a_i N_0}}
\ee
so 
\be
\log \BF_{12}(x) &= \log \frac{b}{a}-\frac{2 (x_2-T)+S}{a N_0}+\frac{2 (x_2-T)+S}{b N_0} \\
&= \log \frac{b}{a}+\frac{(a-b)(2 (x_2-T)+S)}{ab N_0}
\ee

\item \emph{Case 6: $0 < T < T+S < x_2 < x_1$}. In this case the likelihood under $H_i$ is
\be
L(x_2,x_1 \mid \Nt(t)) &= \frac{3}{\Nt(x_2)} \frac{1}{\Nt(x_1)} e^{-2\Lambda(T)-\frac{2 S}{a_i N_0}-2\Lambda(T+S,x_2)-\Lambda(T)-\frac{S}{a_i N_0}-\Lambda(T+S,x_1) }  \\
&= \frac{3}{\Nt(x_2)} \frac{1}{\Nt(x_1)} e^{-3\Lambda(T)-2\Lambda(T+S,x_2)-\Lambda(T+S,x_1)}  e^{-\frac{3 S}{a_i N_0}}
\ee
so 
\be
\log \BF_{12}(x) &= -\frac{3 S}{a N_0}+\frac{3 S}{b N_0} = \frac{3(a-b)S}{ab N_0}
\ee
\end{enumerate} 

\be 
\log \BF_{12}(x) &= \left\{ \begin{array}{lc} 0 & 0 < x_2 < x_1 < T \\
\log \frac{b}{a} + \frac{(a-b)(x_1-T)}{ab N_0} & 0 < 
x_2 < T < x_1 < T+S \\
\frac{(a-b)S}{ab N_0}  & 0 < x_2 < T < T+S < x_1 \\
2 \log \frac{b}{a} + \frac{(a-b)(x_1+2x_2-3T)}{a b N_0} & 0 < T < x_2 < 
x_1  < T+S \\
\log \frac{b}{a} + \frac{(a-b)(2 x_2 -2T+S)}{ab N_0} & 0 < T < x_2 < T+S < x_1 
\\
 \frac{3(a-b)S}{ab N_0} & 0 < T < T+S < x_2 < x_1 
\end{array} \right.
\ee

We go line by line calculating the components of $\P[\vartheta(X) = \vartheta \mid H_1]$. Designate each of the six pieces of the expression by $Q_j$, $j=1,2,\ldots,6$. 

\be
Q_1 = \frac12 \P[X_1<T] &= \frac12 \int_0^T \int_{x_2}^T  \frac3{\Nt(x_2) \Nt(x_1)} e^{-2 \Lambda(x_2)-\Lambda(x_1)} \\
&= \frac12 \int_0^T (e^{-\Lambda(x_2)}-e^{-\Lambda(T)}) 
\frac{3}{\Nt(x_2)} e^{-2 \Lambda(x_2)} dx_2 \\
&= \frac12 \int_0^T \frac{3}{\Nt(x_2)} e^{-3 \Lambda(x_2)} dx_2 - e^{-\Lambda(T)} \int_0^T \frac{3}{\Nt(x_2)} e^{-2 \Lambda(x_2)} dx_2 \\
&=\frac12 \left( 1-e^{-3 \Lambda(T)} - e^{-\Lambda(T)} \frac32 \int_0^T \frac{2}{N_1(x_2)} e^{-2 \Lambda(x_2)} dx_2 \right) \\
&=\frac12 \left( 1-e^{-3 \Lambda(T)} - e^{-\Lambda(T)} \frac32 (1-e^{-2\Lambda(T)}) \right) \\
&=\frac12 +\frac14 e^{-3 \Lambda(T)} - \frac34 e^{-\Lambda(T)} 
\ee
Now define
\be
\delta = \frac{ab N_0}{a-b} \log \frac{a}{b}
\ee
then we have
\be
Q_2 &= \int_0^T \int_T^{T+S} \frac3{\Nt(x_2) \Nt(x_1)} e^{-2\Lambda(x_2)-\Lambda(x_1)} \bone\left\{ \log \frac{b}{a} + 
\frac{(a-b)(x_1-T)}{ab N_0} > 0 \right\} dx_1 dx_2 \\
&= \int_0^T \int_{T+(\delta \wedge S)}^{T+S} \frac3{\Nt(x_2) \Nt(x_1)} e^{-2\Lambda(x_2)-\Lambda(x_1)} dx_1 dx_2 \\
&= (e^{-\Lambda(T)-\frac{\delta \wedge S}{aN_0}} - e^{-\Lambda(T)-\frac{S}{a N_0}}) \frac32 \int_0^T  \frac2{\Nt(x_2)} e^{-2\Lambda(x_2)} dx_2 \\
&= (e^{-\frac{\delta \wedge S}{a N_0}} - e^{-\frac{S}{a N_0}})e^{-\Lambda(T)} \frac32 (1-e^{-2\Lambda(T)})
\ee

For case 3
\be
Q_3 &= \bone\{a>b\} \int_0^T \int_{T+S}^{\infty} \frac{3}{\Nt(x_2) \Nt(x_1)} e^{-2\Lambda(x_2)-\Lambda(x_1)} dx_1 dx_2 \\
&= \bone\{a>b\} e^{-\Lambda(T)-\frac{S}{aN_0}} \frac32 \int_0^T \frac{2}{\Nt(x_2)} e^{-2\Lambda(x_2)} dx_2 \\
&= \frac32 \bone\{a>b\} e^{-\Lambda(T)-\frac{S}{a N_0}} (1-e^{-2 \Lambda(T)}) 
\ee

Case 4
\be
Q_4 &= \int_T^{T+S} \int_{x_2}^{T+S} \frac{3}{\Nt(x_2) \Nt(x_1)} e^{-2\Lambda(x_2)-\Lambda(x_1)} \bone\left\{x_1>T+2(T-x_2+\delta) \right\} dx_1 dx_2 \\
&=  \int_T^{T+S} \int_{x_2}^{T+S} \frac{3}{a^2 N_0^2} e^{-3\Lambda(T)} e^{-\frac{2x_2+x_1-3T}{a N_0}} \bone\left\{x_1>T+2(T-x_2+\delta) \right\} dx_1 dx_2. 
\ee
The inequalities
\be
0&<T<x_2<x_1<T+S, \quad x_1&>T+2(T-x_2+\delta)
\ee
reduce to
\be
\frac{2 \delta}3 &< S < 2\delta, \quad \frac13(3 T + 2 \delta) &< x_1 < S+T, \quad \frac12(3T-x_1+2 \delta) &< x_2 < x_1
\ee
or
\be
S &> 2\delta \\
\intertext{and either}
\frac13(3 T + 2 \delta) &< x_1 < T+2\delta, \quad \frac12(3T-x_1+2 \delta) < x_2 < x_1 \\
\intertext{or}
T + 2 \delta &< x_1 < S+T, \quad T < x_2 < x_1. 
\ee
So then we can express $Q_4$ as
\be
Q_4 &= \begin{cases} 0 & S < \frac23 \delta \\ Q_{41} & \frac23 \delta < S < 2 \delta \\ Q_{42} & S>2\delta \end{cases}
\ee
where
\be
Q_{41} &= \int_{\frac13 (3T+2\delta)}^{T+S} \int_{\frac12(3T-x_1+2\delta)}^{x_1} \frac{3}{a^2 N_0^2} e^{-3\Lambda(T)} e^{-\frac{2x_2+x_1-3T}{a N_0}} dx_2 dx_1 \\
&= \frac12 e^{-3 \Lambda(T)} \left( e^{-\frac{3 S}{a N_0}} - e^{-\frac{2 \delta}{a N_0}} \frac{a N_0 - 3S + 2 \delta}{a N_0} \right)
\ee
and
\be
Q_{42} &= \int_{\frac13 (3T+2\delta)}^{T+2 \delta} \int_{\frac12(3T-x_1+2\delta)}^{x_1} \frac{3}{a^2 N_0^2} e^{-3\Lambda(T)} e^{-\frac{2x_2+x_1-3T}{a N_0}} dx_2 dx_1 \\
&+ \int_{T+2\delta}^{T+S} \int_T^{x_1} \frac{3}{a^2 N_0^2} e^{-3\Lambda(T)} e^{-\frac{2x_2+x_1-3T}{a N_0}} dx_2 dx_1 \\
&= \frac12 e^{-3 \Lambda(T)} \left( e^{-\frac{6 \delta}{a N_0}} + e^{-\frac{2 \delta}{a N_0}} \frac{4 \delta-a N_0}{a N_0} \right) \\
&+ \frac12 e^{-3 \Lambda(T)} \left( e^{-\frac{3 S}{a N_0}}-3e^{-\frac{S}{a N_0}} - e^{-\frac{6 \delta}{a N_0}} + 3 e^{-\frac{2 \delta}{a N_0}} \right) \\
&= \left(\frac{1}{2} e^{-3 \Lambda(T)} \left(e^{-\frac{2 \delta }{a N_0}} \left(\frac{4 \delta }{a N_0}+2\right)+e^{-\frac{3 S}{a N_0}}-3 e^{-\frac{S}{a N_0}}\right)\right)
\ee

And now for case 5
\be
Q_5 &=  \int_T^{T+S} \int_{T+S}^{\infty} f_1(x_1,x_2) \bone\left\{ \log \frac{b}{a}+\frac{(a-b)(2 (x_2-T)+S)}{ab N_0} > 0 \right\} dx_1 dx_2 \\
&=  \int_T^{T+S} \int_{T+S}^{\infty} \frac{3}{\Nt(x_2) \Nt(x_1)} e^{-2\Lambda(x_2)- \Lambda(x_1)} \bone\left\{ x_2 > T + \frac{\delta}2 - \frac{S}2 \right\} dx_1 dx_2 \\
&= \int_{T+\{0 \vee ((\frac{\delta}2 - \frac{S}2) \wedge S)\}}^{T+S} \frac{3}{\Nt(x_2)} e^{-2\Lambda(x_2)} dx_2 \int_{T+S}^{\infty}  \frac{1}{\Nt(x_1)} e^{- \Lambda(x_1)}  dx_1  \\
&= \frac32 e^{-\Lambda(T)-\frac{S}{a N_0}} \int_{T+\{0 \vee ((\frac{\delta}2 - \frac{S}2) \wedge S)\}}^{T+S} \frac{2}{\Nt(x_2)} e^{-2\Lambda(x_2)} dx_2 \\
&= \frac32 e^{-\Lambda(T)-\frac{S}{a N_0}} (e^{-2\Lambda(T+\{0 \vee ((\frac{\delta}2 - \frac{S}2) \wedge S)\})}-e^{-2\Lambda(T+S)}) \\
&= \frac32 e^{-3\Lambda(T)-\frac{S}{a N_0}} (e^{-\frac{2\{0 \vee ((\frac{\delta}2 - \frac{S}2) \wedge S)\}}{aN_0}}-e^{-\frac{2S}{aN_0}})
\ee

Finally case 6
\be
Q_6 &= \bone\{a>b\} \int_{T+S}^{\infty} \int_{x_2}^{\infty} \frac{3}{\Nt(x_2)} \frac{1}{\Nt(x_1)} e^{-2\Lambda(x_2)-\Lambda(x_1)} dx_1 dx_2 \\
&= \bone\{a>b\} \int_{T+S}^{\infty} \frac{3}{\Nt(x_2)} e^{-3\Lambda(x_2)} dx_2 \\
&= \bone\{a>b\} e^{-3\Lambda(T)-\frac{3 S}{aN_0}}
\ee

%

Now we can get the other component fairly easily. We repeat the calculations conditioning on $H_2$
\be
Q_1 &=\frac12 \left( 1+\frac12 e^{-3 \Lambda(T)} - \frac32 
e^{-\Lambda(T)}\right)
\ee

case 2
\be
Q_2 &= \int_0^T \int_T^{T+S} \frac3{\Nt(x_2) \Nt(x_1)} 
e^{-2\Lambda(x_2)-\Lambda(x_1)} \bone\left\{ \log \frac{b}{a} + 
\frac{(a-b)(x_1-T)}{ab N_0} < 0 \right\} dx_1 dx_2 \\
&= \int_0^T \int_{T}^{T+(\delta \wedge S)} \frac3{\Nt(x_2) \Nt(x_1)} 
e^{-2\Lambda(x_2)-\Lambda(x_1)} dx_1 dx_2 \\
&= (e^{-\Lambda(T)} - e^{-\Lambda(T)-\frac{\delta \wedge S}{bN_0}}) \frac32 
\int_0^T  \frac2{\Nt(x_2)} e^{-2\Lambda(x_2)} dx_2 \\
&= (1 - e^{-\frac{\delta \wedge 
S}{bN_0}}) e^{-\Lambda(T)} \frac32 (1-e^{-2\Lambda(T)})
\ee

For case 3
\be
Q_3 &= \bone\{b>a\} \int_0^T \int_{T+S}^{\infty} \frac{3}{\Nt(x_2) \Nt(x_1)} 
e^{-2\Lambda(x_2)-\Lambda(x_1)} dx_1 dx_2 \\
&= 0
\ee

Case 4 
\be
Q_4 &= \int_T^{T+S} \int_{x_2}^{T+S} \frac{3}{\Nt(x_2) \Nt(x_1)} 
e^{-2\Lambda(x_2)-\Lambda(x_1)} \bone\left\{x_1<T+2(T-x_2+\delta) 
\right\} dx_1 dx_2 \\
&=  \int_T^{T+S} \int_{x_2}^{T+S} \frac{3}{b^2 N_0^2} e^{-3\Lambda(T)} 
e^{-\frac{2x_2+x_1-3T}{b N_0}} \bone\left\{x_1<T+2(T-x_2+\delta) 
\right\} dx_1 dx_2.
\ee
The inequalities
\be
0 \le T \le x_2 \le x_1 \le T+S, \quad x_1 \le T+2(T-x_2+\delta), \quad \delta > 0
\ee
reduce to
\be
0 &\le S \le \frac{2 \delta}3, \quad T \le x_1 \le S+T, \quad T \le x_2 \le x_1, \quad \text{or} \\
\frac{2 \delta}3 &\le S \le 2 \delta,\quad  \begin{cases} T < x_1 \le \frac13 (3T+2\delta), & T \le x_2 \le x_1 \quad \text{or} \\ \frac13 (3T+2\delta), \le x_1 \le T+S & T \le x_2 \le \frac12(3T-x_1+2\delta) \end{cases} \\
\intertext{or}
S &> 2 \delta, \quad \begin{cases} T < x_1 < \frac13 (3T+2\delta), & T < x_2 < x_1 \text{or} \\ \frac13 (3T+2\delta) \le x_1 \le T+2\delta, & T < x_2 \le \frac12 (3T-x_1+2\delta) \end{cases}
\ee
so
\be
Q_4 &= \begin{cases} Q_{41} & 0 \le S \le \frac{2 \delta}3 \\ Q_{42} & \frac{2 \delta}3 \le S \le 2 \delta \\ Q_{43} & S > 2 \delta \end{cases},
\ee
where
\be
Q_{41} &= \frac{1}{2} e^{-3\Lambda(T)} \left(e^{-\frac{3 S}{b N_0}}-3 e^{-\frac{S}{b N_0}}+2\right) \\
Q_{42} &= \frac{1}{2} e^{-3 \Lambda(T)} \left(\frac{e^{-\frac{2 \delta }{b N_0}} (b N_0+2 \delta -3 S)}{b N_0}-3 e^{-\frac{S}{b N_0}}+2\right) \\
Q_{43} &= \frac{1}{2bN_0} e^{-3\Lambda(T)} e^{-\frac{2 \delta +S}{b N_0}} \left(b N_0 \left(2 e^{\frac{2 \delta +S}{b N_0}}+e^{\frac{S}{b N_0}}-3 e^{\frac{T}{b N_0}}\right)+e^{\frac{S}{b N_0}} (-4 \delta -3 S+3 T)\right)
\ee

Case 5
\be
Q_5 &=  \int_T^{T+S} \int_{T+S}^{\infty} f_1(x_1,x_2) \bone\left\{ \log 
\frac{b}{a}+\frac{(a-b)(2 (x_2-T)+S)}{ab N_0} < 0 \right\} dx_1 dx_2 \\
&=  \int_T^{T+S} \int_{T+S}^{\infty} \frac{3}{\Nt(x_2) \Nt(x_1)} 
e^{-2\Lambda(x_2)- \Lambda(x_1)} \bone\left\{ x_2 < T + \frac{\delta}2 - 
\frac{S}2 \right\} dx_1 dx_2 \\
&= \int_{T}^{T+\{0 \vee ((\frac{\delta}2 - \frac{S}2) \wedge S)\}} 
\frac{3}{\Nt(x_2)} e^{-2\Lambda(x_2)} dx_2 \int_{T+S}^{\infty}  
\frac{1}{\Nt(x_1)} e^{- \Lambda(x_1)}  dx_1  \\
&= \frac32 e^{-\Lambda(T)-\frac{S}{b N_0}} \int_T^{T+\{0 \vee ((\frac{\delta}2 
- \frac{S}2) \wedge S)\}} \frac{2}{\Nt(x_2)} e^{-2\Lambda(x_2)} dx_2 \\
&= \frac32 e^{-\Lambda(T)-\frac{S}{b N_0}} 
(e^{-2\Lambda(T)}-e^{-2\Lambda(T+\{0 \vee ((\frac{\delta}2 - \frac{S}2) \wedge 
S)\})}) \\
&= \frac32 e^{-3\Lambda(T)-\frac{S}{b N_0}} 
(1-e^{-\frac{2\{0 \vee ((\frac{\delta}2 - \frac{S}2) \wedge 
S)\}}{bN_0}}) 
\ee

Case 6
\be
Q_6 &= \bone\{b>a\} \int_{T+S}^{\infty} \int_{x_2}^{\infty} \frac{3}{\Nt(x_2)} 
\frac{1}{\Nt(x_1)} e^{-2\Lambda(x_2)-\Lambda(x_1)} dx_1 dx_2 \\
&= 0
\ee

\section{Proof of theorem \ref{thm:BayesError2_1Locus}}

\begin{proof}
Define $\Lambda_i(t) = \int_0^t \frac1{N_i(s)} ds$, we then have
\be
\mathbf P[ \vartheta(X)=1 \mid H_1] &= \int_0^{\infty} \bone\left\{ 
\left(\frac1c-1 \right) \Lambda_1(x) > \log 
\frac{1}{c} \right\} \frac{1}{N_1(x)} e^{-\int_0^x \frac{1}{N_1(t)} 
dt} dx \\
&= \mathbf P\left[ X > \Lambda_1^{-1}\left( \frac{c}{1-c} \log \frac1c 
\right) \mid H_1 \right] \\
&= e^{-\Lambda_1\left (\Lambda_1^{-1} \left( \frac{c}{1-c} \log \frac1c 
\right) \right) } = e^{ \frac{c \log c}{1-c} }  \\
&= c^{\frac{c}{1-c}},
\ee
which implicitly assumed that $c < 1$. Similarly
\be
\mathbf P[ \vartheta(X)=2 \mid H_2] &= \int_0^{\infty} \bone\left\{ 
\Lambda_2(x)-\Lambda_1(x) < \log 
\frac{N_1(x)}{N_2(x)} \right\} \frac{1}{N_2(x)} e^{-\int_0^x \frac{1}{N_2(t)} 
dt} dx \\
&= \int_0^{\infty} \bone\left\{ 
\Lambda_2(x)(c-1) > \log 
c \right\} \frac{1}{N_2(x)} e^{-\int_0^x \frac{1}{N_2(t)} 
dt} dx \\
&= \P\left[ X < \Lambda_2^{-1}\left( \frac1{c-1} \log 
c \right) \mid H_2 \right] = 1-e^{-\Lambda_2\left( \Lambda_2^{-1}\left( \frac1{c-1} \log c \right) \right)} \\
&= 1-c^{\frac1{1-c}}
\ee
so then
\be
\P[\vartheta(X)=\vartheta] = \frac12 c^{\frac{c}{1-c}}+\frac12 \left(1-c^{\frac1{1-c}} \right).
\ee
\end{proof}

\section{Proof of theorem \ref{thm:BayesError2_JLoci}}
\begin{proof}
Define $\Lambda_i(t) = \int_0^t \frac1{N_i(s)} ds$ and notice that 
\be
\BF_{12} &= \prod_{j=1}^J \frac{\frac1{N_1(x^j)} e^{-\int_0^{x^j} \frac1{N_1(t)} dt}}{\frac1{N_2(x^j)} e^{-\int_0^{x^j} \frac1{N_2(t)} dt}} = \prod_{j=1}^J \frac{N_2(x^j) e^{-\int_0^{x^j} \frac1{N_1(t)} dt}}{N_1(x^j) e^{-\int_0^{x^j} \frac1{N_2(t)} dt}} \\
&= c^J e^{-\sum_{j=1}^J \Lambda_1(x^j)-\Lambda_2(x^j)} = c^J e^{-\sum_{j=1}^J \Lambda_1(x^j)\left(1-\frac1c \right) } \\
\log \BF_{12} &= J \log c - \left(1-\frac1c \right) \sum_{j=1}^J \Lambda_1(x^j)
\ee
so then
\be
\P \left[ \log \BF_{12} > 0 \mid H_1 \right] &= \P\left[ J \log c > \left(1-\frac1c \right) \sum_{j=1}^J \Lambda_1(X^j) \mid H_1 \right] \\
&= \P\left[ \left(\frac1c-1 \right) \sum_{j=1}^J \Lambda_1(X^j) > J \log \frac1c  \mid H_1 \right] \\
&= \P\left[ \sum_{j=1}^J \Lambda_1(X^j) > J \frac{c}{1-c} \log \frac1c  \mid H_1 \right]. 
\ee
Since
\be
\P[\Lambda_1(X) > s \mid H_1] &= \P[X > \Lambda_1^{-1}(s)] = e^{-s},
\ee
we have
\be
\P \left[ \log \BF_{12} > 0 \mid H_1 \right] &= \P\left[ W > J \frac{c}{1-c} \log \frac1c \right], 
\ee
where $W$ is the sum of $J$ independent unit rate exponentials, so $W \sim \text{Gamma}(J,1)$ and
\be
\P \left[ \log \BF_{12} > 0 \mid H_1 \right] &= 1-\frac1{\Gamma(J)} \gamma\left( J,J \frac{c}{1-c} \log \frac1c \right), 
\ee
where $\gamma(\alpha,\beta)$ is the lower incomplete Gamma function. Similar calculations give us that
\be
\P \left[ \log \BF_{12} < 0 \mid H_2 \right] &= \P\left[ J \log c <  (c-1) \sum_{j=1}^J \Lambda_2(X^j) \mid H_2 \right] \\
&= \P\left[ \sum_{i=1}^J \Lambda_2(X^j) < J \frac1{c-1} \log c  \mid H_2 \right] \\
&= \P\left[ W < J \frac1{c-1} \log c \right] \\
&= \frac1{\Gamma(J)} \gamma\left( J,J \frac{1}{c-1} \log c \right) \\
\ee
giving us
\be
\P[\vartheta(X^J) = \vartheta] &= \frac12 \left( 1-\frac1{\Gamma(J)} \gamma\left( J,J \frac{c}{1-c} \log \frac1c \right) + \frac1{\Gamma(J)} \gamma\left( J,J \frac{1}{c-1} \log c \right) \right),
\ee
as claimed.
\end{proof}
\section{Proof of Theorem \ref{thm:Hellinger}}
We have
\be \label{eq:DensityOnetime}
f_i(x) = \left\{ \begin{array}{lc} \frac{1}{\Nt(x)} e^{-\int_0^x 
\frac{1}{\Nt(t)} dt}  & x<T \\ 
\frac{1}{a_i N_0} e^{-\int_0^T \frac{1}{\Nt(t)} 
dt} e^{-\frac{x-T}{a_i N_0}} & T \le x < T+S \\
\frac{1}{\Nt(x)} e^{-\int_0^T \frac{1}{\Nt(t)} dt} e^{-\frac{S}{a_i N_0}} 
e^{-\int_{T+S}^x \frac{1}{\Nt(t)} dt} & T+S \le x 
\end{array} \right.
\ee
where $f_i(x)$ is the density of a single coalescent time under $H_i$. Define 
\be
\Delta_{12}(x) \equiv (\sqrt{f_1}(x)-\sqrt{f_2}(x))^2.
\ee
So now we calculate
\be
\int (\sqrt{f_1(x)} - \sqrt{f_2(x)})^2 dx &= \int_0^T \Delta_{12}(x) dx + \int_T^{T+S} \Delta_{12}(x) dx + 
\int_{T+S}^{\infty} \Delta_{12}(x) dx
\ee
clearly the first term on the right is zero so
\be
\int \Delta_{12}(x) dx &= \int_T^{T+S}  \Delta_{12}(x) dx + \int_{T+S}^{\infty}  \Delta_{12}(x) dx.
\ee
Observe 
\be
\int_T^{T+S}  \Delta_{12}(x) dx &= \int_{T}^{T+S} \bigg( 
\frac{1}{\sqrt{a N_0}} e^{-\frac12 \int_0^T \frac{1}{\Nt(t)} 
dt} e^{-\frac12 \frac{x-T}{a N_0}} -  \frac{1}{\sqrt{b N_0}} e^{-\frac12 
\int_0^T \frac{1}{\Nt(t)} 
dt} e^{-\frac12 \frac{x-T}{b N_0}}
\bigg)^2 dx \\
&= e^{- \int_0^T \frac{1}{\Nt(t)} dt} \int_{T}^{T+S} \bigg( 
\frac{1}{\sqrt{a N_0}} e^{-\frac12 \frac{x-T}{a N_0}} -  \frac{1}{\sqrt{b 
N_0}} e^{-\frac12 \frac{x-T}{b N_0}}
\bigg)^2 dx,
\ee
then
\be \label{eq:Term2}
e^{\int_0^T \frac{1}{\Nt(t)} dt} \int_T^{T+S} \Delta_{12}(x) dx &= 2 - e^{-\frac{S}{aN_0}} - 
e^{-\frac{S}{bN_0}} - \frac{4 b (1-e^{-\frac{(a+b) S}{2 ab N_0}}) 
\sqrt{a}}{(a+b) \sqrt{b}},
\ee
and now
\be
\int_{T+S}^{\infty} \Delta_{12}(x) dx &= \int_{T+S}^{\infty} 
\bigg( \frac{1}{\sqrt{\Nt(x)}} e^{-\frac12 \int_0^T \frac{1}{\Nt(t)} dt} 
e^{-\frac12 \frac{S}{a N_0}} 
e^{-\frac12 \int_{T+S}^x \frac{1}{\Nt(t)} dt} \\
&- \frac{1}{\sqrt{\Nt(x)}} e^{-\frac12 \int_0^T \frac{1}{\Nt(t)} dt} 
e^{-\frac12 \frac{S}{b N_0}} 
e^{-\frac12 \int_{T+S}^x \frac{1}{\Nt(t)} dt} \bigg)^2 dx \\
&= \int \left(e^{-\frac12 \frac{S}{a N_0}} - e^{-\frac12 \frac{S}{b N_0}} 
\right)^2 \left(\frac{1}{\sqrt{\Nt(x)}} e^{-\frac12 \int_0^T \frac{1}{\Nt(t)} 
dt} 
e^{-\frac12 \int_{T+S}^x \frac{1}{\Nt(t)} dt}\right)^2 dx \\
&= \left(e^{-\frac12 \frac{S}{a N_0}} - e^{-\frac12 \frac{S}{b N_0}} 
\right)^2 e^{- 
\int_0^T \frac{1}{\Nt(t)} 
dt}  \int_{T+S}^{\infty} \frac{1}{\Nt(x)}
e^{- \int_{T+S}^x \frac{1}{\Nt(t)} dt} dx \\
&= \left(e^{-\frac12 \frac{S}{a N_0}} - e^{-\frac12 \frac{S}{b N_0}} 
\right)^2 e^{- 
\int_0^T \frac{1}{\Nt(t)} 
dt}  \left( -e^{- \int_{T+S}^x \frac{1}{\Nt(t)} dt} \bigg|_{T+S}^{\infty} 
\right) \\
&= \left(e^{-\frac12 \frac{S}{a N_0}} - e^{-\frac12 \frac{S}{b N_0}} 
\right)^2 e^{- 
\int_0^T \frac{1}{\Nt(t)} 
dt}  
\ee
and adding this to \eqref{eq:Term2}
\be
H^2(f_1,f_2) = \int \Delta_{12}(x) dx &=  e^{-\int_0^T 
\frac{1}{\Nt(t)}dt} 
\left(1-e^{-\frac{(a+b) S}{2abN_0}} \right) \frac{(a+b-2 \sqrt{ab})}{a+b} \\
&=  e^{-\int_0^T \frac{1}{\Nt(t)}dt} 
\left(1-e^{-\frac{(a+b) S}{2abN_0}} \right) \frac{(\sqrt{a} - 
\sqrt{b})^2}{a+b}, \label{eq:HellingerBoundOneTime}
\ee
which is the same as the last displayed equation on \citet[p 11]{kim2015can}. 

\bibliographystyle{plainnat}
\bibliography{testing-hetero}

\begin{thebibliography}{31}
\providecommand{\natexlab}[1]{#1}
\providecommand{\url}[1]{\texttt{#1}}
\expandafter\ifx\csname urlstyle\endcsname\relax
  \providecommand{\doi}[1]{doi: #1}\else
  \providecommand{\doi}{doi: \begingroup \urlstyle{rm}\Url}\fi

\bibitem[Beerli and Felsenstein(2001)]{Beerli2001}
Peter Beerli and Joseph Felsenstein.
\newblock Maximum likelihood estimation of a migration matrix and effective
  population sizes in $n$ subpopulations by using a coalescent approach.
\newblock \emph{Proceedings of the National Academy of Sciences}, 98\penalty0
  (8):\penalty0 4563--4568, 2001.

\bibitem[Drummond et~al.(2012)Drummond, Suchard, Xie, and
  Rambaut]{drummond2012bayesian}
A.J. Drummond, M.A. Suchard, D.~Xie, and A.~Rambaut.
\newblock Bayesian phylogenetics with {BEAU}ti and the {BEAST} 1.7.
\newblock \emph{Molecular Biology and Evolution}, 29:\penalty0 1969--1973,
  2012.

\bibitem[Felsenstein(2006)]{Felsenstein2006}
Joseph Felsenstein.
\newblock Accuracy of coalescent likelihood estimates: Do we need more sites,
  more sequences, or more loci?
\newblock \emph{Molecular Biology and Evolution}, 23\penalty0 (3):\penalty0
  691--700, 2006.

\bibitem[Felsenstein and Rodrigo(1999)]{joseph_coalescent_1999}
Joseph Felsenstein and Allen~G Rodrigo.
\newblock {C}oalescent {A}pproaches to {HIV} {P}opulation {G}enetics.
\newblock In \emph{The {E}volution of {HIV}}, pages 233--272. {Johns Hopkins
  University} Press, 1999.

\bibitem[Fu et~al.(2016)Fu, Posth, Hajdinjak, Petr, Mallick, Fernandes,
  Furtw{\"a}ngler, Haak, Meyer, Mittnik, Nickel, Peltzer, Rohland, Slon,
  Talamo, Lazaridis, Lipson, Mathieson, Schiffels, Skoglund, Derevianko,
  Drozdov, Slavinsky, Tsybankov, Cremonesi, Mallegni, G{\'e}ly, Vacca, Morales,
  Straus, Neugebauer-Maresch, Teschler-Nicola, Constantin, Moldovan, Benazzi,
  Peresani, Coppola, Lari, Ricci, Ronchitelli, Valentin, Thevenet, Wehrberger,
  Grigorescu, Rougier, Crevecoeur, Flas, Semal, Mannino, Cupillard, Bocherens,
  Conard, Harvati, Moiseyev, Drucker, Svoboda, Richards, Caramelli, Pinhasi,
  Kelso, Patterson, Krause, P{\"a}{\"a}bo, and Reich]{EuroAge}
Qiaomei Fu, Cosimo Posth, Mateja Hajdinjak, Martin Petr, Swapan Mallick, Daniel
  Fernandes, Anja Furtw{\"a}ngler, Wolfgang Haak, Matthias Meyer, Alissa
  Mittnik, Birgit Nickel, Alexander Peltzer, Nadin Rohland, Viviane Slon, Sahra
  Talamo, Iosif Lazaridis, Mark Lipson, Iain Mathieson, Stephan Schiffels,
  Pontus Skoglund, Anatoly~P. Derevianko, Nikolai Drozdov, Vyacheslav
  Slavinsky, Alexander Tsybankov, Renata~Grifoni Cremonesi, Francesco Mallegni,
  Bernard G{\'e}ly, Eligio Vacca, Manuel R.~Gonz{\'a}lez Morales, Lawrence~G.
  Straus, Christine Neugebauer-Maresch, Maria Teschler-Nicola, Silviu
  Constantin, Oana~Teodora Moldovan, Stefano Benazzi, Marco Peresani, Donato
  Coppola, Martina Lari, Stefano Ricci, Annamaria Ronchitelli,
  Fr{\'e}d{\'e}rique Valentin, Corinne Thevenet, Kurt Wehrberger, Dan
  Grigorescu, H{\'e}l{\`e}ne Rougier, Isabelle Crevecoeur, Damien Flas, Patrick
  Semal, Marcello~A. Mannino, Christophe Cupillard, Herv{\'e} Bocherens,
  Nicholas~J. Conard, Katerina Harvati, Vyacheslav Moiseyev, Doroth{\'e}e~G.
  Drucker, Ji{\v r}{\'\i} Svoboda, Michael~P. Richards, David Caramelli, Ron
  Pinhasi, Janet Kelso, Nick Patterson, Johannes Krause, Svante P{\"a}{\"a}bo,
  and David Reich.
\newblock The genetic history of ice age europe.
\newblock \emph{Nature}, 534:\penalty0 200--205, 2016.
\newblock URL \url{http://dx.doi.org/10.1038/nature17993}.

\bibitem[Fu and Li(1993)]{Fu693}
Y~X Fu and W~H Li.
\newblock Statistical tests of neutrality of mutations.
\newblock \emph{Genetics}, 133\penalty0 (3):\penalty0 693--709, 1993.

\bibitem[Gao and Keinan(2016)]{GaoKeinan}
Feng Gao and Alon Keinan.
\newblock Explosive genetic evidence for explosive human population growth.
\newblock \emph{Current Opinion in Genetics and Development}, 41\penalty0
  (Supplement C):\penalty0 130 -- 139, 2016.
\newblock Genetics of human origin.

\bibitem[Gattepaille et~al.(2016)Gattepaille, G{\"u}nther, and
  Jakobsson]{Gattepaillegenetics.115.185058}
Lucie Gattepaille, Torsten G{\"u}nther, and Mattias Jakobsson.
\newblock Inferring past effective population size from distributions of
  coalescent times.
\newblock \emph{Genetics}, 2016.
\newblock \doi{10.1534/genetics.115.185058}.

\bibitem[Griffiths and Tavar{\'e}(1994)]{griffiths_sampling_1994}
R~C Griffiths and S~Tavar{\'e}.
\newblock Sampling theory for neutral alleles in a varying environment.
\newblock \emph{Philosophical Transactions of the Royal Society of London.
  Series B, Biological Sciences}, 344:\penalty0 403--410, June 1994.

\bibitem[Griffiths and Marjoram(1997)]{griffiths1997ancestral}
Robert~C. Griffiths and Paul Marjoram.
\newblock An ancestral recombination graph.
\newblock In Peter Donnelly and Simon Tavar\'{e}, editors, \emph{Progress in
  population genetics and human evolution}, volume~87 of \emph{IMA Volumes in
  Mathematics and Its Applications}, pages 257--270. Springer Verlag, New York,
  1997.

\bibitem[Iles et~al.(2014)Iles, Raghwani, Harrison, Pepin, Djoko, Tamoufe,
  LeBreton, Schneider, Fair, Tshala, Kayembe, Muyembe, Edidi-Basepeo, Wolfe,
  Simmonds, Klenerman, and Pybus]{HCVLast}
James~C Iles, Jayna Raghwani, GL~Abby Harrison, Jacques Pepin, Cyrille~F Djoko,
  Ubald Tamoufe, Matthew LeBreton, Bradley~S Schneider, Joseph~N Fair, Felix~M
  Tshala, Patrick~K Kayembe, Jean~Jacques Muyembe, Samuel Edidi-Basepeo,
  Nathan~D Wolfe, Peter Simmonds, Paul Klenerman, and Oliver~G Pybus.
\newblock Phylogeography and epidemic history of hepatitis {C} virus genotype 4
  in {A}frica.
\newblock \emph{Virology}, 464-465\penalty0 (100):\penalty0 233--243, 09 2014.

\bibitem[Kim et~al.(2015)Kim, Mossel, R{\'a}cz, and Ross]{kim2015can}
Junhyong Kim, Elchanan Mossel, Mikl{\'o}s~Z R{\'a}cz, and Nathan Ross.
\newblock Can one hear the shape of a population history?
\newblock \emph{Theoretical population biology}, 100:\penalty0 26--38, 2015.

\bibitem[Kingman(1982)]{Kingman:1982uj}
John F.~C. Kingman.
\newblock {The coalescent}.
\newblock \emph{Stochastic Processes and Their Applications}, 13\penalty0
  (3):\penalty0 235--248, 1982.

\bibitem[Kuhner et~al.(1995)Kuhner, Yamato, and Felsenstein]{Kuhner:1995vw}
M.K. Kuhner, J.~Yamato, and J.~Felsenstein.
\newblock {Estimating effective population size and mutation rate from sequence
  data using Metropolis-Hastings sampling}.
\newblock \emph{Genetics}, 140\penalty0 (4):\penalty0 1421--1430, 1995.

\bibitem[Li and Durbin(2011)]{Li:2011ez}
Heng Li and Richard Durbin.
\newblock {Inference of human population history from individual whole-genome
  sequences}.
\newblock \emph{Nature}, 475\penalty0 (7357):\penalty0 493--496, July 2011.

\bibitem[Marjoram and Wall(2006)]{MarjoramSMC}
Paul Marjoram and Jeff Wall.
\newblock Fast ``coalescent'' simulation.
\newblock \emph{BMC Genetics}, 7\penalty0 (1), 2006.

\bibitem[McVean and Cardin(2005)]{McVean2005b}
G.~McVean and N.~Cardin.
\newblock {Approximating the coalescent with recombination}.
\newblock \emph{Philos Trans R Soc Lond B Biol Sci}, 360\penalty0
  (1459):\penalty0 1387--1393, Jul 2005.

\bibitem[Myers et~al.(2008)Myers, Fefferman, and Patterson]{MYERS2008342}
Simon Myers, Charles Fefferman, and Nick Patterson.
\newblock Can one learn history from the allelic spectrum?
\newblock \emph{Theoretical Population Biology}, 73\penalty0 (3):\penalty0 342
  -- 348, 2008.

\bibitem[Palacios et~al.(2015)Palacios, Wakeley, and
  Ramachandran]{Palaciosgenetics}
Julia~A. Palacios, John Wakeley, and Sohini Ramachandran.
\newblock Bayesian nonparametric inference of population size changes from
  sequential genealogies.
\newblock \emph{Genetics}, 201\penalty0 (1):\penalty0 281--304, 2015.

\bibitem[Pluzhnikov and Donnelly(1996)]{PluzhnikovDonnelly}
Anna Pluzhnikov and Peter Donnelly.
\newblock Optimal sequencing strategies for surveying molecular genetic
  diversity.
\newblock \emph{Genetics}, 144:\penalty0 1247--1262, 1996.

\bibitem[Sainudiin et~al.(2011)Sainudiin, Thornton, Harlow, Booth, Stillman,
  Yoshida, Griffiths, McVean, and Donnelly]{Sainudiin2011}
Raazesh Sainudiin, Kevin Thornton, Jennifer Harlow, James Booth, Michael
  Stillman, Ruriko Yoshida, Robert Griffiths, Gil McVean, and Peter Donnelly.
\newblock Experiments with the site frequency spectrum.
\newblock \emph{Bulletin of Mathematical Biology}, 73\penalty0 (4):\penalty0
  829--872, 2011.

\bibitem[Schiffels and Durbin(2014)]{MSMC}
Stephan Schiffels and Richard Durbin.
\newblock Inferring human population size and separation history from multiple
  genome sequences.
\newblock \emph{Nature Genetics}, 46\penalty0 (8):\penalty0 919--925, 2014.

\bibitem[Shapiro et~al.(2004)Shapiro, Drummond, Rambaut, Wilson, Matheus, Sher,
  Pybus, Gilbert, Barnes, Binladen, Willerslev, Hansen, Baryshnikov, Burns,
  Davydov, Driver, Froese, Harington, Keddie, and Kosintsev]{shapiro_rise_2004}
Beth Shapiro, Alexei~J. Drummond, Andrew Rambaut, Michael~C. Wilson, Paul~E.
  Matheus, Andrei~V. Sher, Oliver~G. Pybus, M.~Thomas~P. Gilbert, Ian Barnes,
  Jonas Binladen, Eske Willerslev, Anders~J. Hansen, Gennady~F. Baryshnikov,
  James~A. Burns, Sergei Davydov, Jonathan~C. Driver, Duane~G. Froese,
  C.~Richard Harington, Grant Keddie, and Pavel Kosintsev.
\newblock Rise and fall of the {B}eringian steppe bison.
\newblock \emph{Science}, 306\penalty0 (5701):\penalty0 1561--1565, 2004.

\bibitem[Sheehan et~al.(2013)Sheehan, Harris, and Song]{diCal}
Sara Sheehan, Kelley Harris, and Yun~S. Song.
\newblock Estimating variable effective population sizes from multiple genomes:
  A sequentially {M}arkov conditional sampling distribution approach.
\newblock \emph{Genetics}, 194\penalty0 (3):\penalty0 647--662, 2013.

\bibitem[Slatkin and Hudson(1991)]{Slatkin:1991wx}
M.~Slatkin and R.R. Hudson.
\newblock {Pairwise comparisons of mitochondrial DNA sequences in stable and
  exponentially growing populations}.
\newblock \emph{Genetics}, 129\penalty0 (2):\penalty0 555--562, 1991.

\bibitem[Stephens and Donnelly(2000)]{StephensDonnelly2000}
Matthew Stephens and Peter Donnelly.
\newblock Inference in molecular population genetics.
\newblock \emph{Journal of the Royal Statistical Society: Series B (Statistical
  Methodology)}, 62\penalty0 (4):\penalty0 605--635, 2000.

\bibitem[Terhorst and Song(2015)]{Terhorst23062015}
Jonathan Terhorst and Yun~S. Song.
\newblock Fundamental limits on the accuracy of demographic inference based on
  the sample frequency spectrum.
\newblock \emph{Proceedings of the National Academy of Sciences}, 112\penalty0
  (25):\penalty0 7677--7682, 2015.

\bibitem[Terhorst et~al.(2017)Terhorst, Kamm, and Song]{TerhorstSFS}
Jonathan Terhorst, John~A. Kamm, and Yun~S. Song.
\newblock Robust and scalable inference of population history from hundreds of
  unphased whole genomes.
\newblock \emph{Nature Genetics}, 49:\penalty0 303--309, 2017.

\bibitem[Tong et~al.(2015)Tong, Shi, LiuDi, Qian, Liang, Bo, Liu, Ren, Fan, Ni,
  Sun, Jin, Teng, Li, Kargbo, Dafae, Kanu, Chen, Lan, Jiang, Luo, Lu, Zhang,
  Yang, Hu, Cao, Deng, Su, Sun, Liu, Wang, Wang, Bu, Guo, Zhang, Nie, Bai, Sun,
  An, Xu, Zhang, Huang, Mi, Yu, Yao, Feng, Xia, Zheng, Yang, Lu, Jiang, Kargbo,
  He, Gao, Cao, and in~{S}ierra {L}eone]{EbolaBeast}
Yi-Gang Tong, Wei-Feng Shi, LiuDi, Jun Qian, Long Liang, Xiao-Chen Bo, Jun Liu,
  Hong-Guang Ren, Hang Fan, Ming Ni, Yang Sun, Yuan Jin, Yue Teng, Zhen Li,
  David Kargbo, Foday Dafae, Alex Kanu, Cheng-Chao Chen, Zhi-Heng Lan, Hui
  Jiang, Yang Luo, Hui-Jun Lu, Xiao-Guang Zhang, Fan Yang, Yi~Hu, Yu-Xi Cao,
  Yong-Qiang Deng, Hao-Xiang Su, Yu~Sun, Wen-Sen Liu, Zhuang Wang, Cheng-Yu
  Wang, Zhao-Yang Bu, Zhen-Dong Guo, Liu-Bo Zhang, Wei-Min Nie, Chang-Qing Bai,
  Chun-Hua Sun, Xiao-Ping An, Pei-Song Xu, Xiang-Li-Lan Zhang, Yong Huang,
  Zhi-Qiang Mi, Dong Yu, Hong-Wu Yao, Yong Feng, Zhi-Ping Xia, Xue-Xing Zheng,
  Song-Tao Yang, Bing Lu, Jia-Fu Jiang, Brima Kargbo, Fu-Chu He, George~F. Gao,
  Wu-Chun Cao, and The China Mobile Laboratory Testing~Team in~{S}ierra
  {L}eone.
\newblock Genetic diversity and evolutionary dynamics of {E}bola virus in
  {S}ierra {L}hateone.
\newblock \emph{Nature}, 524\penalty0 (7563):\penalty0 93--96, 2015.

\bibitem[Wakeley(2008)]{wakeley_coalescent_2008}
John Wakeley.
\newblock \emph{Coalescent Theory: An Introduction}.
\newblock Roberts \& Company Publishers, June 2008.

\bibitem[Wakeley and Sargsyan(2009)]{Wakeley10112008}
John Wakeley and Ori Sargsyan.
\newblock Extensions of the coalescent effective population size.
\newblock \emph{Genetics}, 181:\penalty0 341--345, 2009.

\end{thebibliography}

\end{document}